%% file: Main_Inventiones.tex
\documentclass[12pt]{article}

\usepackage{amsfonts,pstricks}
\usepackage{amsmath, amsthm, amssymb,  mathrsfs, epsfig}
\usepackage{hyperref}
\usepackage{graphicx}
\usepackage{tikz}
\usepackage{caption}
\usepackage{authblk}
\usepackage{bm}
\usepackage{tensor}
\usetikzlibrary{decorations.markings,intersections}
\usepackage{enumitem}
\hypersetup{linktocpage = true}
\usepackage{tikz-cd}

\input{Macros} 

\begin{document}

\title{On Two Invariants of Three Manifolds from Hopf Algebras}
\renewcommand{\thefootnote}{\fnsymbol{footnote}} 
\renewcommand{\thefootnote}{\arabic{footnote}} 
\author{Liang Chang%
       \thanks{Email: \texttt{changliang996@nankai.edu.cn}}}
\affil{Chern Institute of Mathematics,\\ Nankai University, Tianjin 300071, China}

\author{Shawn X. Cui%
       \thanks{Email: \texttt{xingshan@stanford.edu}}}
\affil{Stanford Institute for Theoretical Physics,\\ Stanford University, Stanford, California 94305}
\maketitle

\begin{abstract}
We prove a $20$-year-old conjecture concerning two quantum invariants of three manifolds that are constructed from finite dimensional Hopf algebras, namely, the Kuperberg invariant and the Hennings-Kauffman-Radford invariant. The two invariants can be viewed as a non-semisimple generalization of the Turaev-Viro-Barrett-Westbury $(\text{TVBW})$ invariant  and the Witten-Reshetikhin-Turaev $(\text{WRT})$ invariant, respectively. By a classical result relating $\text{TVBW}$ and  $\text{WRT}$, it follows that the Kuperberg invariant for a semisimple Hopf algebra is equal to the Hennings-Kauffman-Radford invariant for the Drinfeld double of the Hopf algebra. However, whether the relation holds for non-semisimple Hopf algebras has remained open, partly because the introduction of framings in this case makes the Kuperberg invariant significantly more complicated to handle. We give an affirmative answer to this question. An important ingredient in the proof involves using a special Heegaard diagram  in which one family of circles gives the surgery link of the three manifold represented by the Heegaard diagram. 
\end{abstract}

\section{Introduction}
\label{sec:introduction}
Since the discovery of the Jones polynomial \cite{jones1985polynomial} and the formulation of a topological quantum field theory ($\TQFT$)\cite{witten1988topological} \cite{atiyah1988topological} in the $1980s$, there haven been fascinating interactions between low dimensional topology and quantum physics. Many quantum invariants of $3$-manifolds have been constructed, which deeply connects together different areas of research such as knot theory, tensor categories, quantum groups, Chern-Simons theory, conformal field theory, etc. 
Quantum invariant generally refers to the partition function of a $\TQFT$, or less rigorously, to any invariant that is defined as a state-sum model. In dimension three, tensor categories and Hopf algebras are the main sources for quantum invariants. For instance, the Turaev-Viro-Barrett-Westbury invariant $Z_{\text{TVBW}}$ \cite{turaev1992state} \cite{barrett1996invariants} and the $\wrt$ invariant $Z_{\text{WRT}}$ \cite{reshetikhin1991invariants} are based on spherical fusion categories and modular categories, respectively. Both invariants can be extended to a $\TQFT$ and the latter is believed to be a mathematical realization of Witten-Chern-Simons theory. These invariants are particularly important in topology as they distinguish certain homotopy equivalent $3$-manifolds \cite{sokolov1997lens}.

Two fundamental invariants that are constructed from finite dimensional Hopf algebras in the early $1990s$ are the Kuperberg invariant $Z_{\text{Kup}}$ \cite{kuperberg1991involutory} \cite{kuperberg1996non} and the $\hkr$ invariant $Z_{\text{HKR}}$ \cite{Hennings} \cite{Kauffman-Radford}. On one hand, $Z_{\text{Kup}}$ is defined for any finite dimensional Hopf algebra and is an invariant of framed oriented closed $3$-manifolds. If the Hopf algebra is semisimple, then $Z_{\text{Kup}}$ does not depend on framings and hence becomes an invariant of closed oriented $3$-manifolds. On the other hand, the $Z_{\text{HKR}}$ invariant, initially defined by Hennings and later reformulated by Kauffman and Radford, is an invariant of closed oriented $3$-manifolds, but can be naturally refined to also include a $2$-framing (similar to $Z_{\text{WRT}}$). Moreover, $Z_{\text{HKR}}$ requires the Hopf algebra to be ribbon in addition to some non-degeneracy conditions (see Section \ref{subsec:hkr}). 

The $Z_{\text{HKR}}$ invariant has been extensively studied in the literature. In \cite{lyubashenko1995invariants} \cite{lyubashenko1995modular}, Lyubashenko produced an invariant from certain monoidal categories (not necessarily semisimple) which generalized both $Z_{\text{HKR}}$ and $Z_{\text{WRT}}$. The relation between $Z_{\text{HKR}}$ and $Z_{\text{WRT}}$ for semisimple Hopf algebras and certain quantum groups were explored in \cite{kerler1996genealogy} \cite{chen2009relation} \cite{chen2012three} \cite{doser2017relationship} \cite{kerler2003homology}.  $\TQFT$ properties of  $Z_{\text{HKR}}$ were given in \cite{kerler1998connectivity} \cite{chen2017integrality} \cite{de2017renormalized}.  Murakami combined ideas from $Z_{\text{HKR}}$ and $Z_{\text{WRT}}$ to define a generalized Kashaev invariant of links in $3$-manifolds and proposed a version of volume conjecture for this invariant \cite{murakami2017generalized}.

It has been a long-standing conjecture that $Z_{\text{Kup}}$ from a Hopf algebra $H$ is equal to $Z_{\text{HKR}}$ from the Drinfeld double $D(H)$ of $H$, namely, for any closed oriented $3$-manifold $X$, 
\begin{align}
\label{equ:Kup=HKR semisimple}
Z_{\text{Kup}}(X;H) = Z_{\text{HKR}}(X;D(H)).
\end{align}
The relation was speculated in \cite{kuperberg1996non} and stated explicitly (and more generally for Lyubashenko invariant) in \cite{kerler1996genealogy}\footnote{The issue of framings was not mentioned in both of these references, but we will address it below.}. Since then, there have been many partial results along this direction. Barrett and Westbury proved \cite{barrett1995equality} that for  semisimple $H$,
\begin{equation}
\label{equ:Kup=TVBW}
Z_{\text{Kup}}(X;H) = Z_{\text{TVBW}}(X; \text{Rep}(H)).
\end{equation} 
Similarly Kerler \cite{kerler1996genealogy} proved that for semisimple and modular $H$,
\begin{equation}
\label{equ:HKR=WRT}
Z_{\text{HKR}}(X;H) = Z_{\text{WRT}}(X;\text{Rep}(H)).
\end{equation} 
In this sense, $Z_{\text{Kup}}$ and $Z_{\text{HKR}}$ can be considered as non-semisimple generalizations of $Z_{\text{TVBW}}$ and $Z_{\text{WRT}}$, respectively. If $\CatC$ is a spherical fusion category, then the $\Drinfeld$ double $D(\CatC)$ of $\CatC$ is a modular category. Turaev and Virelizier  \cite{turaev2010two} proved 
\begin{equation}
\label{equ:TVBW=WRT}
Z_{\text{TVBW}}(X;\CatC) = Z_{\text{WRT}}(X;D(\CatC)),
\end{equation}
which generalizes the well-known result for the case of $\CatC$ modular \cite{walker1991witten} \cite{turaev1992quantum} \cite{Roberts}
\begin{equation}
\label{equ:TVBW=WRT2}
Z_{\text{TVBW}}(X;\CatC) = Z_{\text{WRT}}(X \# \overline{X};\CatC).
\end{equation}
Equations \ref{equ:Kup=TVBW} \ref{equ:HKR=WRT} \ref{equ:TVBW=WRT} together imply the conjecture in Equation \ref{equ:Kup=HKR semisimple} for semisimple Hopf algebras. A direct proof of the conjecture in this case was also given by Sequin in his thesis \cite{sequin2012comparing}. However, whether Equation \ref{equ:Kup=HKR semisimple} holds for non-semisimple Hopf algebras has remained to be a somewhat $20$-year-old open problem.  Another consequence implied from the categorical counterpart and also conjectured in \cite{kerler1996genealogy} is that when $H$ itself is ribbon and semisimple, we have
\begin{equation}
\label{equ:Kup=HKR2 semisimple}
Z_{\text{Kup}}(X;H) = Z_{\text{HKR}}(X \# \overline{X};H),
\end{equation}
and again this has been verified directly in \cite{Chang2015}. In the current paper, we aim to give a proof of (a suitable variation) of both Equation \ref{equ:Kup=HKR semisimple} and \ref{equ:Kup=HKR2 semisimple} for non-semisimple Hopf algebras. Explicitly, we prove the following two theorems.
 
\begin{theorem}
\label{thm:main11}
Let $H$ be a finite dimensional $\db$ Hopf algebra and $X$ be a closed oriented $3$-manifold, then there exist a framing $b$ and a $2$-framing $\phi$ of $X$ such that,
\begin{align}
\label{equ:main11}
\Kup{X,b;H} = \HKR{X,\phi; D(H)}.
\end{align}
\end{theorem}

\begin{theorem}
\label{thm:main22}
Let $H$ be a finite dimensional factorizable ribbon Hopf algebra and $X$ be a closed oriented $3$-manifold, then there exists a framing $b$ of $X$ such that
\begin{align}
\Kup{X,b;H}=\HKR{X\#\overline{X};H}.
\end{align}
\end{theorem}

One feature of the paper is an extensive use of tensor diagrams in computing both $Z_{\text{Kup}}$ and $Z_{\text{HKR}}$. In fact, both invariants can be defined by tensor diagrams alone. This implies that the results in the current paper not only hold for Hopf algebras in the category of vector spaces, but also hold for Hopf super-algebras or Hopf objects in a monoidal category which sufficiently resembles the category of finite dimensional vector spaces. For the sake of simplicity, we restrict the discussions on ordinary Hopf algebras.

These two theorems reveal a connection between Hopf algebras and $3$-manifolds, which is expected to be extended as some type of duality in the category level. In one direction, Hopf algebras yield topological invariants of $3$-manifolds; in the other direction, we can study Hopf algebras using topology. When the $3$-manifold is fixed, $Z_{\text{Kup}}$ and $Z_{\text{HKR}}$ may provide algebraic invariants for Hopf algebras. Two Hopf algebras are said to be gauge equivalent if their representation categories are equivalent as tensor categories or equivalently, they are connected by the twisting of some 2-cocycle. One family of gauge invariants are the Frobenius-Schur indicators \cite{kashina2012trace} \cite{kashina2006higher}, which have important applications to the representation theory and coincide with $Z_{\text{Kup}}$ for lens space \cite{Chang-Wang}. It is speculated that $Z_{\text{Kup}}$ provides more general gauge invariants for any finite dimensional Hopf algebras. By a recent result on gauge dependence of  $Z_{\text{HKR}}$ (\cite{chen2017integrality}), Theorem \ref{thm:main22} implies that $Z_{\text{Kup}}$ is a gauge invariant for ribbon Hopf algebras. More detailed discussions will appear in a subsequent paper.   

One issue that is not solved here is whether the $2$-framing on the $\text{RHS}$ of Equation \ref{equ:main11} is the same as the one induced by the framing on the $\text{LHS}$. Since a change of $2$-framing by one unit changes the $Z_{\text{HKR}}$ by a root of unity, this issue is not relevant up to roots of unity. Another question is whether Equation \ref{equ:main11} still holds for all framings $b$ and the corresponding $\phi$ induced from $b$. We leave it as a future direction.

The rest of the paper is organized as follows. In Section \ref{sec:Hopf} we give a review and set up the conventions on Hopf algebra. Some Lemmas on Hopf algebras will be proved for use later. Section \ref{sec:invariants} recalls the definition of the invariants $Z_{\text{Kup}}$ and $Z_{\text{HKR}}$. In particular, we refine the latter to include $2$-framings. Section \ref{sec:main1} and \ref{sec:main2} are devoted to the proof of our main results, Theorem \ref{thm:main11} and Theorem \ref{thm:main22}, respectively.

\section{Hopf Algebras}
\label{sec:Hopf}
In this section we give a minimal review on Hopf algebras and prove a few lemmas. For a detailed treatment of Hopf algebras, see, for instance, \cite{kuperberg1996non} \cite{radford1994trace} \cite{radford2012hopf}, etc. Formulas in Hopf algebras are illustrated either by tensor diagrams or algebraic expressions. It is straight forward to convert one notation into the other. A novelty in this section is to represent the structure maps in the $\Drinfeld$ double by tensor diagrams from the original Hopf algebra, which turns out convenient to manipulate relations in the double and useful later in comparing different invariants of $3$-manifolds.   Throughout the context, Let $H = H(M, i,\Delta, \epsilon, S)$ be a finite dimensional Hopf algebra over $\Complex$, where the symbols inside the parenthesis denote the multiplication, unit, comultiplication, counit, and antipode, respectively.

\subsection{Tensor Networks}
\label{subsec:tensor}
\begin{figure}
\centering
\includegraphics[scale=1]{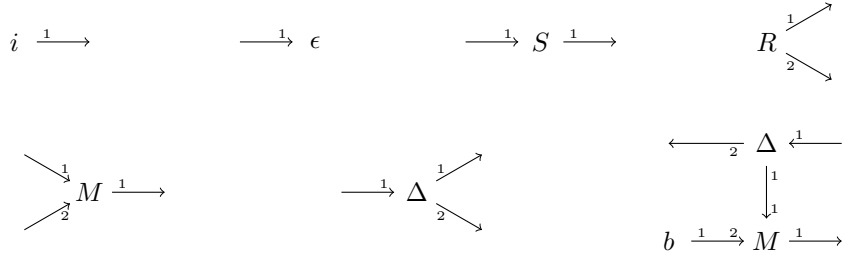}
\caption{Examples of tensor diagrams}\label{fig:tensor}
\end{figure}

Tensor networks have wide applications in physics and quantum information. For a review of tensor networks, see \cite{orus2014practical} \cite{cui2016quantum}, etc.  In \cite{kuperberg1991involutory} \cite{kuperberg1996non}, tensor networks are used as a convenient tool to represent and manipulate operations in Hopf algebras. Let $V$ be a finite dimensional vector space and $V^*$ be its dual. A {\it tensor diagram} in $V$ is a pair $(\CatG, \CatT = \{\CatT_v\})$ where,
\begin{itemize}
\item $\CatG$ is a directed graph such that at each vertex $v$, there is a local ordering on the set of incoming legs (i.e., edges) and a local ordering on the set of outgoing legs by $\{1,\cdots, i_v\}$ and $\{1, \cdots, o_v\}$, respectively; 
\item for each vertex $v$, $\CatT_v \in V^{1} \otimes \cdots \otimes V^{i_v} \otimes V_{1} \otimes \cdots \otimes V_{o_v}$, where each $V^i$ is a copy of $V^*$ associated with the $i$-th incoming leg and each $V_j$ is a copy of $V$ associated with the $j$-th outgoing leg. In this case, $\CatT_v$ is called an $(i_v, o_v)$ tensor. 
\end{itemize}
Choose a basis $\{v_1, \cdots, v_k\}$ of $V$ and a dual basis $\{v^1, \cdots, v^k\}$ of $V^*$, then an $(m,n)$ tensor $\CatT$ can be written as
\begin{align}
\CatT = \sum \CatT_{i_1, \cdots, i_m}^{j_1, \cdots, j_n}\, v^{i_1} \otimes \cdots \otimes v^{i_m} \otimes v_{j_1} \otimes \cdots \otimes v_{j_n}.
\end{align}
See Figure \ref{fig:tensor} for examples of tensor diagrams on the plane. In these diagrams, vertices are replaced by the labels of the corresponding tensors. Around a vertex, a number is placed beside each leg to represent the local ordering. An $(m,n)$ tensor can be equivalently viewed as a linear map from $V^{\otimes m}$ to $V^{\otimes n}$. From this perspective, a $(0,1)$ tensor is a vector, a $(1,0)$ tensor is a co-vector, a $(1,1)$ tensor is a linear map from $V$ to $V$, etc. Let $(\CatG, \CatT)$ be a tensor diagram. Assume there are $i_{\CatG} + o_{\CatG}$ dangling legs, $i_{\CatG}$ of them incoming and $o_{\CatG}$ outgoing (with a local ordering of each set), then a contraction of the tensors along all {\it internal} legs results in an $(i_{\CatG}, o_{\CatG})$ tensor, which we call the evaluation of $(\CatG, \CatT)$. By abuse of language, we do not distinguish a tensor diagram with its evaluation.   

Now we make an important convention to simplify drawing tensor diagrams. At each vertex of a tensor diagram, we always group the incoming legs and the outgoing legs. Unless noted otherwise, incoming legs are enumerated {\it counter clockwise} and outgoing legs {\it clockwise}. This uniquely determines a local ordering if both types of legs are present: 
\begin{center}
\includegraphics[scale=1]{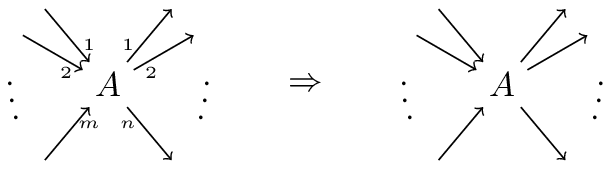}
\end{center}
If there is only one type of legs and the tensor is neither a $(1,0)$ tensor nor a $(0,1)$ tensor, we mark the leg labeled by $1$ explicitly to avoid ambiguities:
\begin{center}
\includegraphics[scale=1]{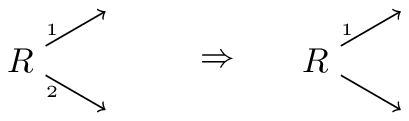}
\end{center}

If $V = H$ is a Hopf algebra, the tensor diagrams of the structure maps are represented by those with the corresponding labels in Figure \ref{fig:tensor}. Relations of between these maps can also be illustrated in tensor diagrams. For instance, the equations $(\Delta \otimes id)\circ \Delta = (id \otimes \Delta)\circ \Delta$ and $\Delta \circ M = (M \otimes M)\circ(id \otimes P \otimes id)\circ (\Delta \otimes \Delta)$ are represented by:
\begin{center}
\includegraphics[scale=1]{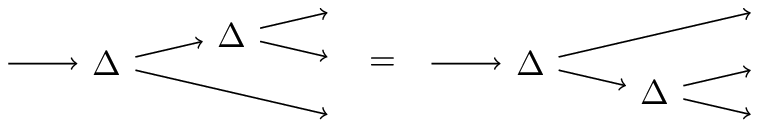}
\qquad
\includegraphics[scale=1]{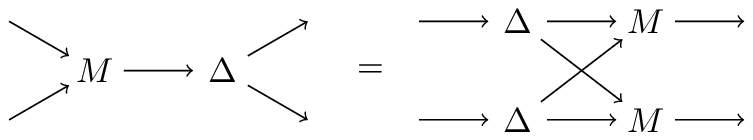}
\end{center}
where $P: H \otimes H \longrightarrow H \otimes H$ is the swap map. For $n \geq 1$, denote the tensor diagrams for the maps $(\Delta \otimes id^{\otimes (n-2)}) \circ \cdots \circ (\Delta \otimes id) \circ \Delta$ and $M \circ (M \otimes id) \circ \cdots \circ (M \otimes id^{\otimes (n-2)})$ by:
\begin{equation}
\label{equ:general_Delta_M}
\vcenter{\hbox{\includegraphics[scale=1]{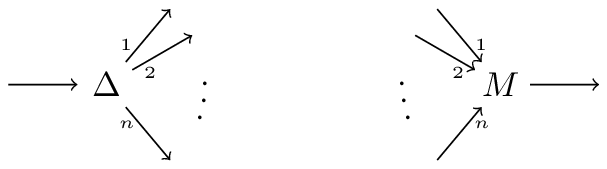}}}
\end{equation}


An $(m,n)$ tensor $\CatT$ in $V$ can also be viewed as an $(n,m)$ tensor $\CatT^*$ in $V^*$ by
\begin{align}
{\CatT^*}^{i_1, \cdots, i_m}_{j_1, \cdots, j_n}:= \CatT_{i_1, \cdots, i_m}^{j_1, \cdots, j_n}.
\end{align}
If $\CatT$ is interpreted as a map from $V^{\otimes m}$ to $V^{\otimes n}$, then $\CatT^*$ is the dual map of $\CatT$. For instance, if $V = H$ is a Hopf algebra, then $\Delta^*$ is a $(2,1)$ tensor representing the multiplication in $V^*$ and for $f, f' \in V^*$, $\Delta^*(f \otimes f')$ is given by: 
\begin{equation}
\vcenter{\hbox{\includegraphics[scale=1]{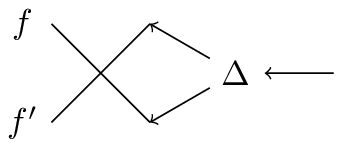}}}
\end{equation}
Note that there is a swap of the two outgoing legs in the above diagram because of our convention for the implicit ordering of the incoming/outgoing legs. The dual notion of tensors will be used in Section \ref{subsubsec:double} when dealing with the quantum double of Hopf algebras.

\subsection{Integrals in Hopf Algebras}
\label{subsubsec:integrals}
A left (resp. right) integral of $H$ is an element $e_L \in H$ (resp. $e_R \in H$) such that $x e_L = \epsilon(x)e_L$ (resp. $e_R x = \epsilon(x)e_R$) for any $x \in H$. Left and right integrals of $H^*$ are denoted by $\mu_L$ and $\mu_R$, respectively. The defining equations of $e_L, \,e_R, \,\mu_L, $ and $\mu_R$ \footnote{In \cite{kuperberg1996non} they are called left cointegral, right cointegral, left integral and right integral, respectively, of $H$.} in terms of tensor diagrams are given by:
\begin{equation}
\label{equ:integral_def}
\vcenter{\hbox{\includegraphics[scale=1]{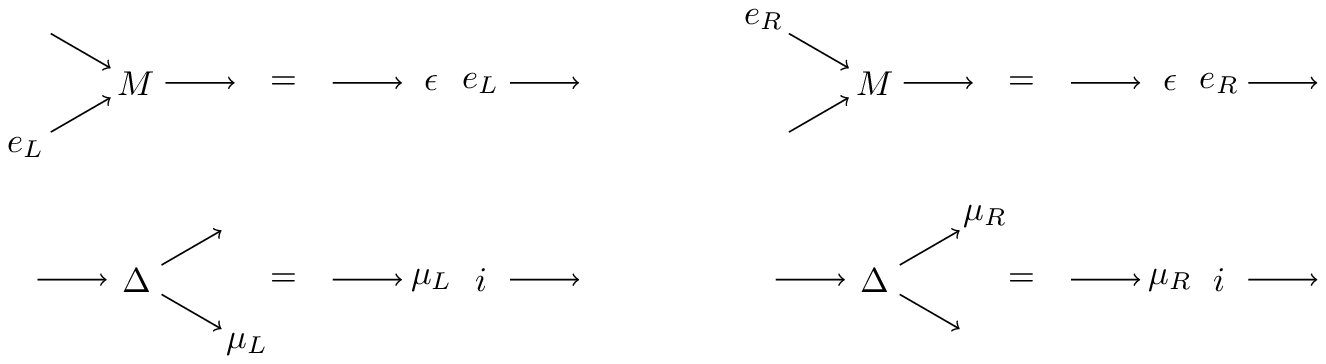}}}
\end{equation}

The space of left integrals and the space of right integrals are both one dimensional. Choose right integrals $e_R \in H,\, \mu_R \in H^*$ such that $\mu_R(e_R) = 1$. Define the distinguished group-like elements $a \in H, \alpha \in H^*$ by,
\begin{equation}
\label{equ:a_alpha_def}
\vcenter{\hbox{\includegraphics[scale=1]{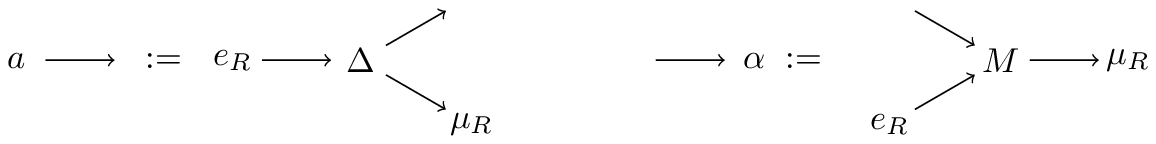}}}
\end{equation}
and for $n \in \Integer$ define $\mu_{n - \frac{1}{2}} \in H^*, \, e_{n - \frac{1}{2}} \in H$ by
\begin{equation}
\label{equ:general_integral_def}
\vcenter{\hbox{\includegraphics[scale=1]{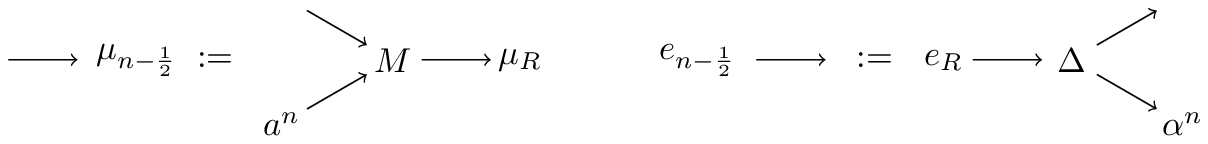}}}
\end{equation}
Then $\mu_R = \mu_{-\frac{1}{2}},\, e_R = e_{-\frac{1}{2}}$ are right integrals and $\mu_L:= \mu_{\frac{1}{2}} \in H^*, \, e_L:= e_{\frac{1}{2}} \in H$ are left integrals. Set $q:= \alpha(a)$. It follows that $q$ is a root of unity and we have $\mu_R(e_R) = \mu_R(e_L) = \mu_L(e_R) = 1$ and $\mu_L(e_L) = q^{-1}$. Moreover, $\mu_L \circ S = \mu_R,$ $\mu_R \circ S = q \mu_L,$ $S(e_L) = e_R,$ $S(e_R) = qe_L$. Thus $S^2$ has eigenvalue $q$ on all integrals of $H$ and $H^*$. The relations between integrals and the distinguished group-like elements are given as follows:
\begin{equation}
\label{equ:a_alpha_integral_rel}
\vcenter{\hbox{\includegraphics[scale=1]{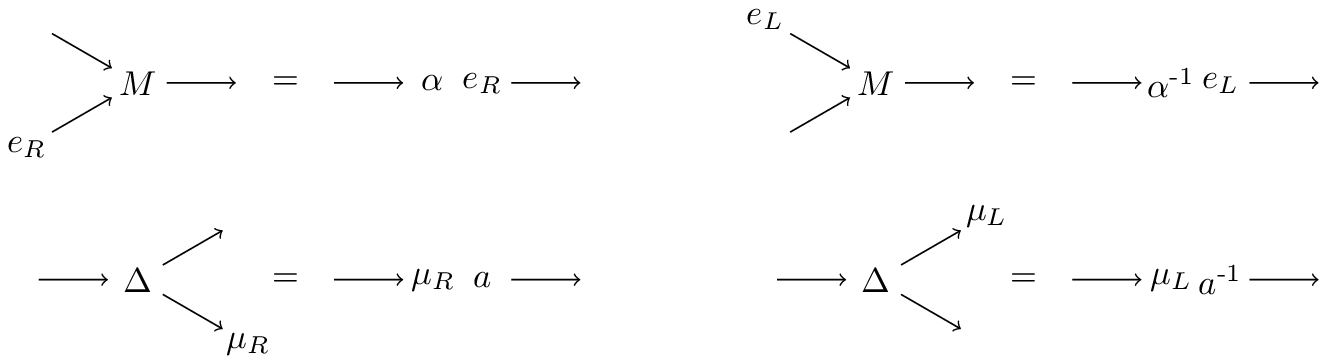}}}
\end{equation}

Note that here $a$ and $\alpha$ correspond to $g$ and $\alpha^{-1}$, respectively, in \cite{radford1994trace} \cite{radford2012hopf}.
The well-known Radford formula for $S^4$ can be expressed as
\begin{equation}
\label{equ:S4_def}
\vcenter{\hbox{\includegraphics[scale=1]{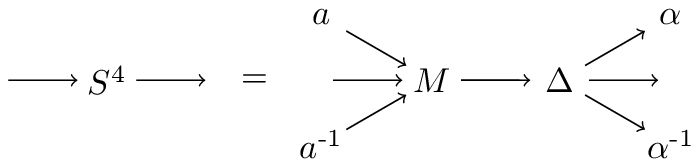}}}
\end{equation}
Also define
\begin{equation}
\label{equ:T_def}
\vcenter{\hbox{\includegraphics[scale=1]{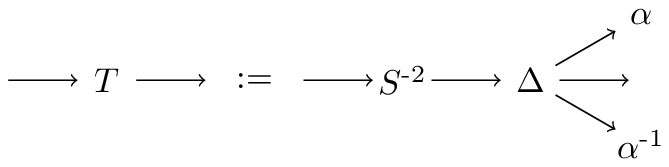}}}
\end{equation}
Then $T$ is an automorphism of $H$ as a Hopf algebra, i.e., $T$ commutes with all structure maps of $H$.
\begin{lemma}
\label{lem:T and S}
For any $n \in \Integer$,
\begin{itemize}
\item $S^2$ has eigenvalue $q$ on $e_{n-\frac{1}{2}}$ and $\mu_{n-\frac{1}{2}}$, namely, $S^2(e_{n-\frac{1}{2}}) = q e_{n-\frac{1}{2}}$, $\mu_{n-\frac{1}{2}} \circ S^2 = q \mu_{n-\frac{1}{2}}$.
\item $T$ fixes $a, \, \alpha, \,e_{n-\frac{1}{2}},\, \mu_{n-\frac{1}{2}}$, namely, $T(a) =a$, $\alpha \circ T = \alpha$, $T(e_{n-\frac{1}{2}}) = e_{n-\frac{1}{2}}$, $\mu_{n-\frac{1}{2}} \circ T = \mu_{n-\frac{1}{2}}.$
\end{itemize}
\begin{proof}
The first part follows directly from the calculation:
\begin{center}
\includegraphics[scale=1]{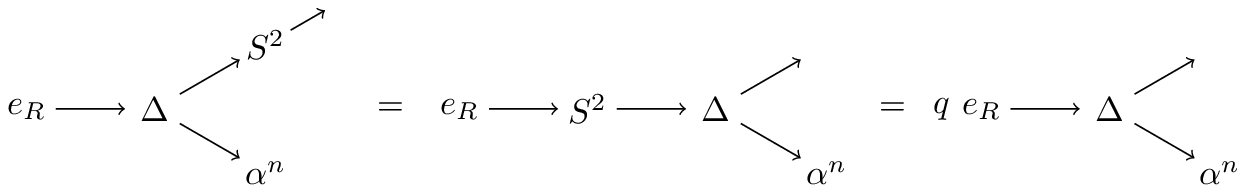}
\end{center}

For the second part, $\mu_R \circ (S^2T)$ is computed as follows:
\begin{center}
\includegraphics[scale=1]{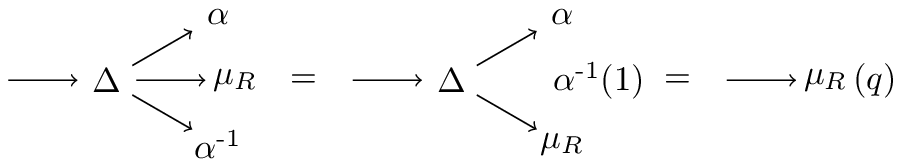}
\end{center}
where the first equality is by definition of $\mu_R$ and the second equality is by Equation \ref{equ:a_alpha_integral_rel}. Hence $\mu_R \circ (S^2T) = q \mu_R$. By the first part, we have $\mu_R \circ T = \mu_R$. 

By using the Radford formua in Equation \ref{equ:S4_def}, $S^2T^{-1}$ can be expressed as
\begin{center}
\includegraphics[scale=1]{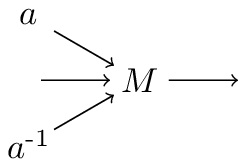}
\end{center}
A similar calculation as above shows that $S^2T^{-1} (e_R) = q e_R$, and thus $T(e_R) = e_R$. That $T$ also fixes $e_{n-\frac{1}{2}} $ and $ \mu_{n-\frac{1}{2}}$ follows immediately.
\end{proof}
\end{lemma}

The Hopf algebra $H$ is called {\it balanced} if $T = id$, and {\it unimodular} if left integrals of $H$ are also right integrals. The latter is equivalent to the condition that $\alpha = \epsilon$. If $H$ is unimodular, then $q = 1$, $e_L = e_R \in Z(H)$, and for any $x,y \in H$,  we have
\begin{align}
\mu_R \circ S^2 (x) = \mu_R(x), \qquad \mu_R(xy) = \mu_R\left(S^2(y)x\right) = \mu_R\left(yS^{-2}(x)\right). 
\end{align} 

\subsection{Ribbon Hopf Algebras}
\label{subsubsec:ribbon}
A {\it $\qt$} Hopf algebra is a pair $(H,R)$, where $H$ is a finite dimensional Hopf algebra, $R \in H \otimes H$, called the $R$-matrix, is an invertible element, and for any $x \in H$,
\begin{align}
R \Delta(x) = \Dop(x)R, \quad (\Delta \otimes id)(R) = R_{13}R_{23}, \quad (id \otimes \Delta)(R) = R_{13}R_{12}. 
\end{align}
If $(H, R)$ is $\qt$, then
{\small
\begin{align}
R^{-1} = (S \otimes id)(R) = (id \otimes S^{-1})(R), \quad R = (S \otimes S)(R), \quad (\epsilon \otimes id)(R) = (id \otimes \epsilon)(R) = 1.
\end{align} }
Let $u:= M \circ (S \otimes id )(R_{21}) \in H$ be the {\it $\Drinfeld$ element}. Then $u$ is invertible and $S^2(x) = uxu^{-1}$ for any $x \in H$. Moreover, $S(u)u = uS(u) \in Z(H)$, and if $H$ is unimodular, then $ uS(u)^{-1} = a$ is the distinguished group-like element. Set $Q = R_{21}R$ and define the $\Drinfeld$ map $f_{Q}$ by
\begin{align*}
f_{Q}: H^* \longrightarrow H, \quad p \longmapsto (p \otimes id)Q.
\end{align*}
Then $f_{Q}(\alpha^{-1}) = 1 = f_{Q}(\epsilon)$. 

The pair $(H,R)$ is called {\it factorizable} if $(H,R)$ is $\qt$ and $f_Q$ is a linear isomorphism. Thus factorizable Hopf algebras are unimodular, with the distinguished group element given by $uS(u)^{-1}$. Let $(H,R)$ be factorizable and $\mu_R$ be a right integral, then $f_{Q}(\mu_R)$ is a (two-sided) integral of $H$ and one can choose $\mu_R$ such that $\mu_R \circ f_{Q}(\mu_R) = 1$.   

A ribbon Hopf algebra is a triple $(H, R, v)$ where $(H,R)$ is a $\qt$ Hopf algebra and $v \in Z(H)$, called the ribbon element, satisfies the following equation:
\begin{align}
v^2 = uS(u), \quad S(v) = v, \quad \epsilon(v) = 1, \quad \Delta(v) = (v \otimes v) Q^{-1} .
\end{align}
Since $u$ is invertible, so is $v$. Let $G:= u v^{-1}$. Then $G$ is a group-like element and $G^2 = u^2v^{-2} = uS(u)^{-1}$. 


\subsection{The Quantum Double of Hopf Algebras}
\label{subsubsec:double}
\begin{figure}
\centering
\includegraphics[scale=1]{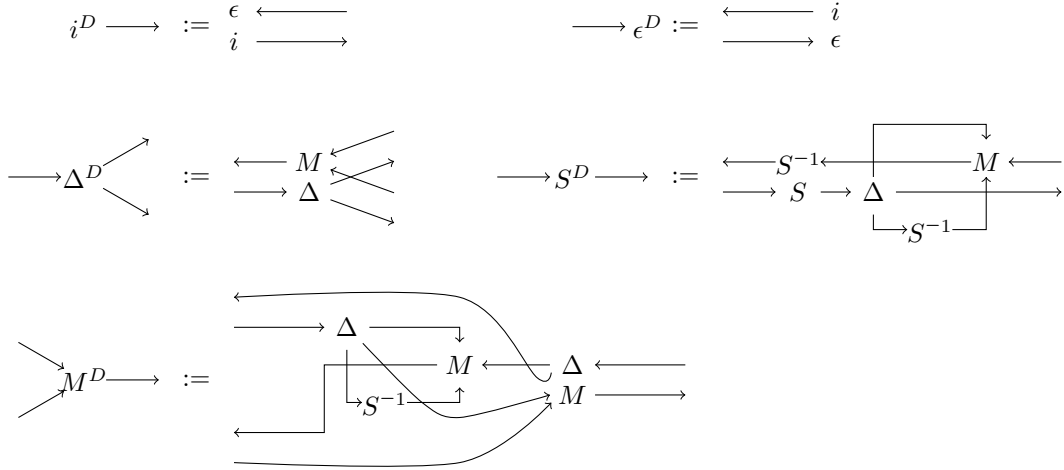}
\caption{Definition of Hopf algebra structures in $D(H)$}\label{fig:double}
\end{figure}

\begin{figure}
\centering
\includegraphics[scale=1]{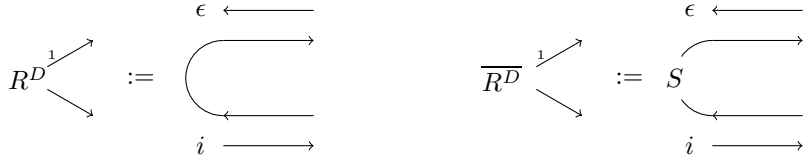}
\caption{The $R$-matrix in $D(H)$, where $\overline{R^D}$ means $(R^D)^{-1}$.}\label{fig:double_RandRinverse}
\end{figure}

Introduced in \cite{drinfeld1986quantum}, the quantum double (or $\Drinfeld$ double) $D(H) = H^{*\cop} \otimes H$ of a Hopf algebra $H$ is a factorizable $\qt$ (and thus unimodular) Hopf algebra. Instead of writing down algebraically the Hopf algebra structures in $D(H)$, we describe them with tensor diagrams consisting of tensors in $H$, which will be used later in Section \ref{sec:main1} to describe the $Z_{\text{HKR}}$ invariant from a quantum double. Labels for operations in the double will be endowed with a superscript \lq$D$'. For instance, $\Delta^D$ means the comultiplication in $D(H)$.  The vector $f \otimes v \in D(H)$ and  covector $v \otimes f \in D(H)^*$ are represented respectively by
\begin{equation}
\vcenter{\hbox{\includegraphics[scale = 1]{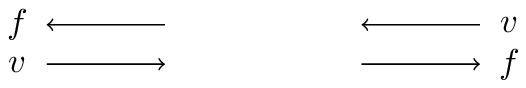}}}
\end{equation}
That is, we use a pair of oppositely directed arrows to represent a copy of $D(H)$ with the arrow on the top corresponding to $H^{*\cop}$ and the one on the bottom to $H$. The definition of the Hopf algebra structures in $D(H)$ using tensor diagrams are given in Figure \ref{fig:double}. Keep in mind that for a tensor with both incoming and outgoing legs, the incoming legs are listed in {\it counter-clockwise} order while the outgoing legs {\it clockwise}. The $R$-matrix is given in Figure \ref{fig:double_RandRinverse}.

One advantage of using tensor diagrams is that it provides a direct visualization on how structures in the double are constructed from those in the original Hopf algebra. It is also convenient for deriving equations. Of course, one can always obtain the algebraic expressions from the diagrams. For instance, for $f \otimes v, f' \otimes v' \in D(H)$, from Figure \ref{fig:double} we see the the multiplication is given by:
\begin{align}
M^D\left((f \otimes v) \otimes (f '\otimes v')\right) =f (v_{(1)} \rightharpoonup f' \leftharpoonup  S^{-1}(v_{(3)}) ) \otimes v_{(2)} v'
\end{align}

With notations from Section \ref{subsubsec:integrals} , let $\mu_R^D = e_L \otimes \mu_R,$ $\mu_L^D = q^{-1} e_R \otimes \mu_L, $ $e^D = q \mu_L \otimes q_R$, or in tensor diagrams,
\begin{equation}
\label{equ:double_integrals}
\vcenter{\hbox{\includegraphics[scale=1]{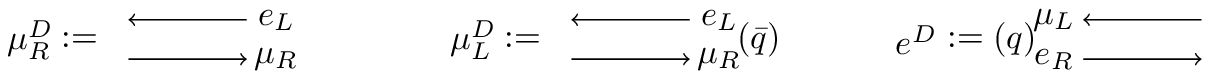}}}
\end{equation}
Then $\mu_R^D$ and $\mu_L^D$ are a right integral and left integral of $D(H)^*$, respectively, and $e^D$ is a two-sided integral of $D(H)$. Moreover, $\mu_R^D(e^D) = \mu_L^D(e^D) = 1$. The distinguished group-like element is given by $a^D = \alpha^{-1} \otimes a$ which can be checked as follows: 
\begin{equation}
\label{equ:double_computing_a}
\vcenter{\hbox{\includegraphics[scale=1]{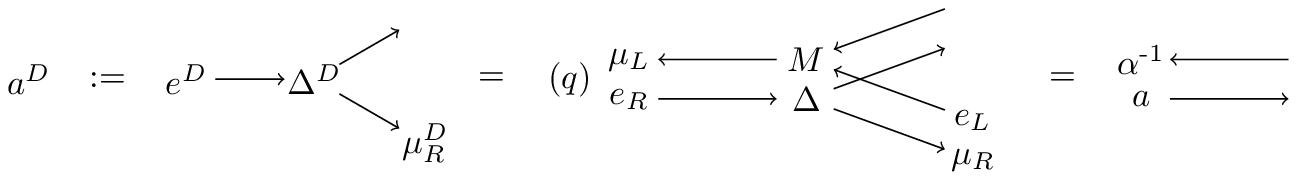}}}
\end{equation}
where the first equality above is by definition and the third equality is from Equation \ref{equ:a_alpha_integral_rel}.

\begin{lemma}
\label{lem:computing_M_tilde}
Let $f \otimes x \in D(H)$ and $n \in \Integer$, then
\begin{align}
\mu_R^D\left((f \otimes x)(a^D)^{-n}\right) = q^n\, \mu_{-n-\frac{1}{2}}(x) f (e_{n+ \frac{1}{2}}). 
\end{align} 
\begin{proof}
This proof is illustrated in Figure \ref{fig:computing_M_tilde}. The first equality is due to the fact that $\alpha$ is an algebra morphism. The second equality uses the definition of $T$ and $e_L$. The third equality follows from Lemma \ref{lem:T and S}.
\begin{figure}
\centering
\includegraphics[scale=1]{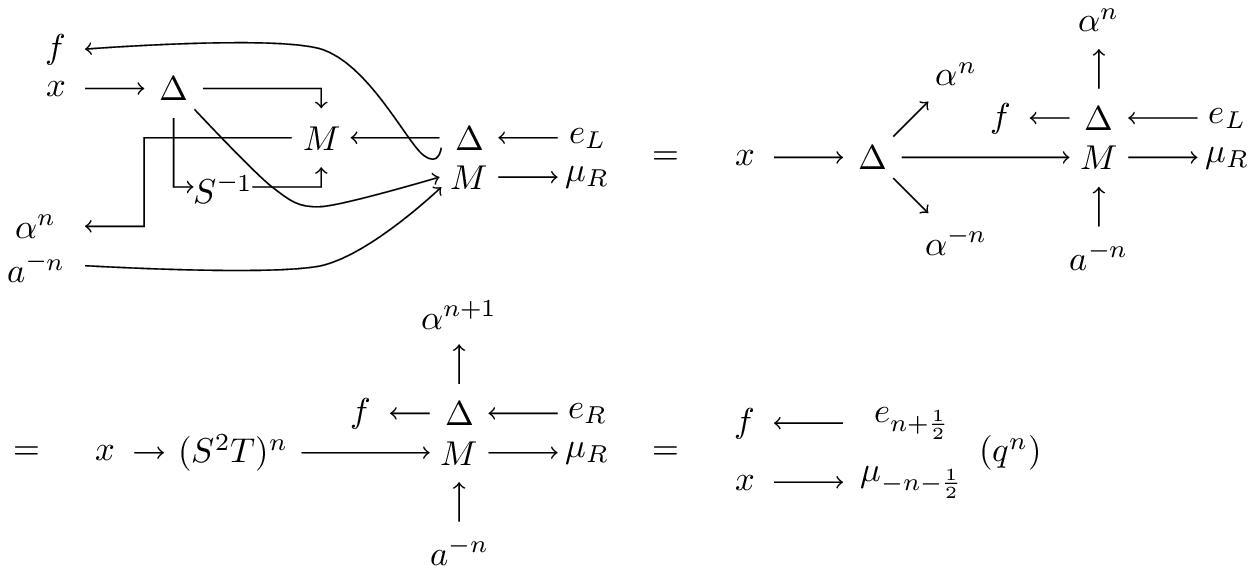}
\caption{Computing $\mu_R^D\left((f \otimes x)(a^D)^{-n}\right)$}\label{fig:computing_M_tilde}
\end{figure}
\end{proof}
\end{lemma}

In general, $D(H)$ may not have ribbon elements. By \cite{kauffman1993necessary}, $D(H)$ is ribbon if and only there exist group-like elements $b \in H, \beta \in H^*$ such that $b^2 = a,$ $\beta^2 = \alpha,$ and,
\begin{equation}
\vcenter{\hbox{\includegraphics[scale=1]{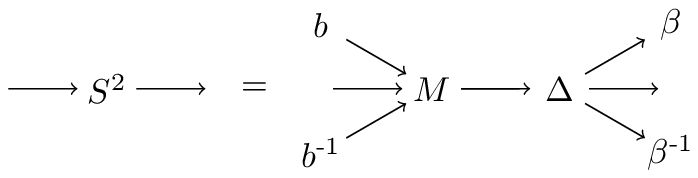}}}
\end{equation}
In \cite{chen2017integrality}, this condition is called {\it $\db$}, but this is not to be confused with the {\it balanced} condition defined in Section \ref{subsubsec:integrals}. It is direct to see that $\tau:= \beta(b)$ is a fourth root of $q$. The corresponding ribbon element of $D(H)$ is given by:
\begin{equation}
\label{equ:double_v}
\vcenter{\hbox{\includegraphics[scale=1]{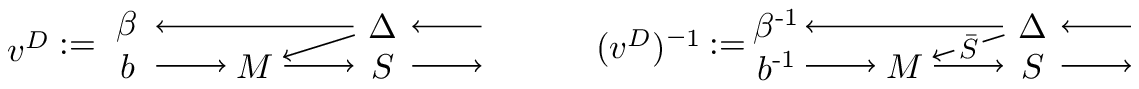}}}
\end{equation}
where $\bar{S}$ in the above diagram means $S^{-1}$.

By direct calculations, $G^D = \beta^{-1} \otimes b$, $\mu_R^D(v^D) = \tau^{-5},$ $\mu_R^{D}\left((v^D)^{-1}\right) = \tau$.

\section{Invariants from Hopf Algebras}
\label{sec:invariants}
In this section we review and make some clarifications on the definitions of the Kuperberg invariant and the $\hkr$ invariant. 

\subsection{Kuperberg Invariant $Z_{\text{Kup}}$}
\label{subsec:Kuperberg}
 The Kuperberg invariant is defined for closed framed oriented $3$-manifolds from a finite dimensional Hopf algebra \cite{kuperberg1996non}. If the Hopf algebra is semi-simple, then the invariant becomes independent of the framings, and is reduced to the invariant of closed oriented $3$-manifolds in \cite{kuperberg1991involutory}. 

We first recall the definitions of combings, framings, and their representations on Heegaard diagrams. Let $X$ be a closed oriented $3$-manifold endowed with a Riemannian metric. A {\it combing} of $X$ is a unit-norm vector field considered up to homotopy, and a {\it framing} of $X$ consists of three orthonormal vector fields consistent with the orientation, again considered up to homotopy. Since the tangent bundle of $X$ is trivial, the set of combings (resp. framings) correspond to homotopy classes of maps from $X$ to $\Sphere^2$ (resp. $\SO(3)$), although the correspondence is in general not canonical. 
Let $R = (\Sigma_g, \alpha, \beta)$ be a Heegaard diagram of $X$ where $\Sigma_g$ is a closed oriented surface of genus $g$, and $\alpha$ and $\beta$ are the collection of lower circles and upper circles, respectively. We only consider minimal Heegaard diagrams. That is, $\alpha$ and $\beta$ each contains exactly $g$ circles. In the following, $R$ and $\Sigma_g$ will be used interchangeably when no confusion arises. Different diagrams of $X$ are related by circle slide, stabilization, and isotopy. Let $\vec{n}$ be the unit normal vector field of $\Sigma_g$ in $X$ pointing from the lower handlebody to the upper handlebody. By convention, the orientation on $\Sigma_g$ and $\vec{n}$ form the orientation on $X$. Any vector field on $X$ can be orthogonally projected along $\vec{n}$ to a tangent vector field on $\Sigma_g$, which could have singularities. The converse problem of extending a vector field on $\Sigma_g$ with certain properties to one on $X$ is studied in \cite{kuperberg1996non}.

\begin{figure}
\centering
\includegraphics[scale=0.7]{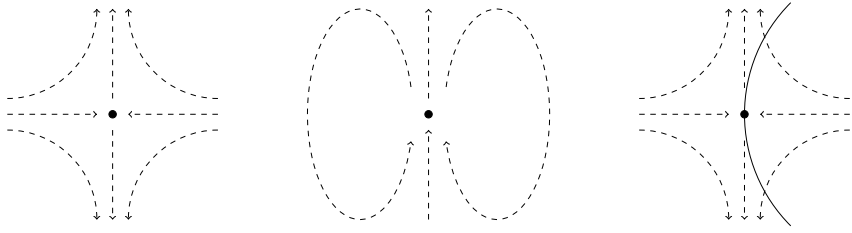}
\caption{(Left) Singularity of index $-1$; (Middle) Singularity of index $2$; (Right) Local picture of a circle (solid curve) near the singularity.}\label{fig:singularity}
\end{figure}

According to \cite{kuperberg1996non}, any combing of $X$ can be represented by a combing of $\Sigma_g$, which, by definition, is a vector field on $\Sigma_g$ with $2g$ singularities of index $-1$, one on each circle, and one more singularity of index $2$ disjoint from all circles. Moreover, each singularity of index $-1$ is distinct from all crossings of the circles, and the two out-pointing vectors should be tangent to the circle. See Figure \ref{fig:singularity} for the local geometry of singularities and the circle near the singularity on it. Any combing $b$ of $\Sigma_g$ can be extended to a combing $\tilde{b}$ of $X$ whose projection to $\Sigma_g$ is the same as $b$, and moreover, one can choose $\tilde{b}$ in such a way that it coincides with $b$ on $\Sigma_g$ away from a small neighborhood of singularities, and at the singularity on a lower (resp. upper) circle $\tilde{b}$ is opposite (resp. parallel) to $\vec{n}$.

\begin{figure}
\centering
\includegraphics[scale=1]{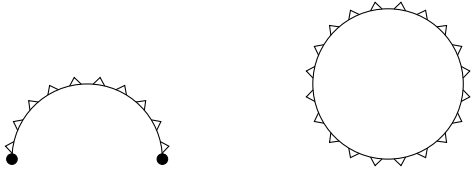}
\caption{An illustration of twist fronts.}\label{fig:fronts1}
\end{figure}

A framing of $X$ is determined by two orthonormal combings $(\tilde{b}_1, \tilde{b}_2)$ since the third one can be inferred from the first two and the orientation. By the previous argument, we can represent the framing as two orthogonal combings $(b_1, b_2)$ on $\Sigma_g$. For reasons that will become clear below, we represent $b_2$ in a different but equivalent form to a combing. Let $\Sigma_g^*$ be the punctured surface of $\Sigma_g$ with all singularities of $b_1$ removed. Then $(b_1, \vec{n}, b_1 \times \vec{n})$ forms an orthogonal frame on $\Sigma_g$ where $b_1 \times \vec{n}$ is the vector orthogonal to both $b_1$ and $\vec{n}$ such that the triple $(b_1, \vec{n}, b_1 \times \vec{n})$ matches the orientation of $X$. Since $b_2$ is orthogonal to $b_1$, $b_2$ lies in the plane spanned by $\vec{n}$ and $b_1 \times \vec{n}$. Then we can define a map $f: \Sigma_g^* \rightarrow \Sphere^1$ by sending $x$ to $f(x) = (f_1(x),f_2(x))$ such that,
\begin{align*}
b_2(x) = f_1(x) \vec{n} + f_2(x) b_1 \times \vec{n}.
\end{align*}
By perturbing $b_2$ in general position, one can assume $(1,0)$ is a regular value of $f$ and hence $f^{-1}(1,0)$ is a $1$-manifold. Namely, the set of points at which $b_2$ is parallel to $\vec{n}$ is a $1$-manifold, where each connected component is either a simple closed curve or an open curve approaching to some singularities in both directions. We also attach small triangles (See figure \ref{fig:fronts1}) on one side of the curves to indicate the direction in which $b_2$ is rotating about $b_1$ by the right-hand rule. More specifically, $f$ takes values in the first quadrant at points which are close to the curve and are located on the side of the curve with triangles. The curves with small triangles attached are called {\it twist fronts}. Twist fronts determine $b_2$ on $\Sigma_g$. Given a collection of twist fronts indicating $b_2$, the following condition needs to be satisfied in order to extend $b_2$ to a combing on $X$ orthogonal to $\tilde{b}_1$.

Arbitrarily orient all circles and consider the frame $(\vec{n}, b_1, \vec{n} \times b_1)$ on $\Sigma_g^*$. For each lower or upper circle $c$, the tangent vector field $c'$  lies in the plane spanned by $b_1$ and $\vec{n} \times b_1$. Define $\theta_{c}$ to be the total counter-clockwise rotation, in unit of $1 = 360^{\circ}$, of $c'$ relative to $b_1$ in the direction of $c$ \footnote{Strictly speaking, the rotation of $c'$ around $b_1$ at the singularity does not make sense since $b_1$ vanishes. Then $\theta_{c}$ is actually defined as the limit $\lim_{x \to *, y \to *} \theta_{c[x,y]}$, where \lq$*$' is the singularity on $c$, $x$ (resp. $y$) is a point of $c$ near the singularity in the forward (resp. backward) direction, and $c[x,y]$ is the subarc of $c$ from $x$ to $y$. Similar situation applies to the definition of $\theta_c(p)$ for a point $p$ on $c$ to be introduced below}. Note that near the singularity $c'$ is parallel to $b_1$ in the forward direction and anti-parallel to $b_1$ in the backward direction, thus $\theta_{c}$ is always a proper half integer. A crossing of a circle with a twist front is called positive if the circle travels from the non-triangle side to the triangle side, and is negative otherwise. Define $\phi_{c}$ to be the number of signed crossings of $c$ with twists fronts, with the crossings at the singularity counted half as much. Then the condition $b_1$ and $b_2$ need to satisfy is:
\begin{align}
\theta_{c} = 
\begin{cases}
\phi_{c}  &, \,c \text{\, is a lower circle } \\
-\phi_{c} &, \,c \text{\, is an upper circle }
\end{cases}
\end{align}

\begin{figure}
\centering
\includegraphics[scale=1]{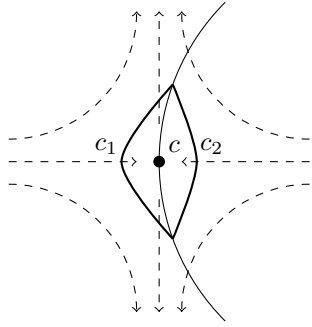}
\caption{Perturbation of a circle $c$ off its singularity.}\label{fig:perturb}
\end{figure}

\begin{remark}
A more intrinsic way to define $\phi_c$ is to use total rotations similar to the definition of $\theta_c$. There are two ways to perturb the circle $c$ off its singularity. See Figure \ref{fig:perturb}. Let $c_1$ and $c_2$ denote the circles resulting from the two perturbations, hence they are contained in $\Sigma_g^*$. Consider the orthogonal frame $(b_1, \vec{n}, b_1 \times \vec{n})$ on $\Sigma_g^*$. Note that $b_2$ lies in the plane spanned by $\vec{n}$ and $b_1 \times \vec{n}$. Define $\phi_{c_i}$ to be the total counter-clockwise rotation of $b_2$ relative to $\vec{n}$ along the curve $c_i$. Then one can check that $\phi_{c_i}$ equals the number of signed crossings of $c_i$ with twist fronts, and the previously defined $\phi_c$ is $ \frac{\phi_{c_1} + \phi_{c_2}}{2}$.  
\end{remark}

Let $p$ be a point on a $c$. Define $\theta_c(p)$ to be the counter-clockwise rotation of $c'$ relative to $b_1$ going along the circle from the singularity to $p$, and define $\phi_c(p)$ to be the number of signed crossings of $c$ with twist fronts from a point near the singularity in the forward direction to $p$. Arrange the diagram so that lower circles intersect upper circles orthogonally. If $p$ is the point of crossing of the lower circle $l$ with the upper circle $u$, let
\begin{align}
\theta(p) := 2(\theta_l(p) - \theta_u(p)) + \frac{1}{2},\quad \phi(p) := \phi_l(p) - \phi_u(p).
\end{align}
It can be shown that $\theta(p)$ is always an integer. Actually, $\theta(p)$ is even if and only if $l$ and $u$ form a positive basis of the tangent space at $p$. We note that in the original definition of $\theta(p)$ in \cite{kuperberg1996non}, the last term is $-\frac{1}{2}$ instead of $\frac{1}{2}$, but we will stick to the current convention as only with this convention, the invariant to be defined will reduce to the one introduced in \cite{kuperberg1991involutory} when the Hopf algebra is semi-simple.

We are ready to define the Kuperberg invariant. Let $H$ be a finite dimensional Hopf algebra. We will use notations from Section \ref{sec:Hopf}. Choose a right integral $\mu_{R}$ and a right co-integral $e_{R}$ so that $\mu_R(e_R) = 1$, and recall the definitions of $\mu_n, e_n$ for $n$ a half integer. Let $X$ be a closed orientated $3$-manifold with a framing $b = (b_1,b_2)$ given on a Heegaard diagram $R = (\Sigma_g, \alpha, \beta)$, where $\alpha = \{\alpha_1, \cdots, \alpha_g\}$ and $\beta = \{\beta_1, \cdots, \beta_g\}$ are lower and upper circles, respectively. Orient all circle arbitrarily, and call the singularity on each circle the basepoint. The definition of the Kuperberg invariant is best illustrated using tensors and tensor contractions. We also given an alternative way to interpreted it afterwards. 

\begin{figure}
\centering
\includegraphics[scale=1]{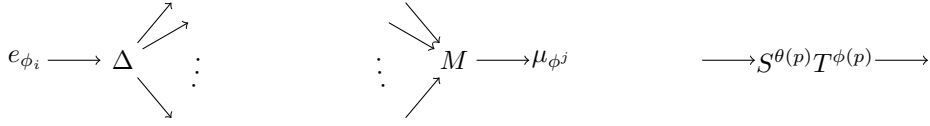}
\caption{(Left) the $\Delta$ tensor assigned to $\alpha_i$; (Middle) the $M$ tensor assigned to $\beta_j$; (Right) the $ST$ tensor assigned to a crossing $p$.}\label{fig:three_tensors}
\end{figure}

For each lower circle $\alpha_i$, let $\phi_i  = \phi_{\alpha_i} (= \theta_{\alpha_i})$ and assign the tensor in Figure \ref{fig:three_tensors}(Left) to $\alpha_i$, one leg for each crossing on $\alpha_i$ counted from the basepoint along its orientation. Similarly for each upper circle $\beta_j$, let $\phi^j = \phi_{\beta_j} (= -\theta_{\beta_j})$ and assign the tensor in Figure \ref{fig:three_tensors}(Middle) to $\beta_j$. For each crossing $p$, insert the tensor shown in Figure \ref{fig:three_tensors}(Right) to connect the two legs, one from the tensor of the lower circle and one from the tensor of the upper circle. Then one obtains a tensor network consisting of the three families of tensors from Figure \ref{fig:three_tensors} without free legs. The Kuperberg invariant $\Kup{X,b;H}$ is then defined to be the contraction of this tensor network.

A more \lq algebraic' but also more lengthy way to define the invariant is as follows. Enumerate the crossings by $p_1, p_2, \cdots, p_m$. Let
\begin{align}
H_{\alpha} = \bigotimes\limits_{i=1}^{g} H(\alpha_i), \quad H_{\beta} = \bigotimes\limits_{i=1}^{g} H(\beta_i), \quad H_{c} = \bigotimes\limits_{i=1}^{m} H(p_i),
\end{align}
Where each $H(\cdot)$ is a copy of $H$.
For each lower circle $\alpha_i$, let $p_{i_1}, \cdots, p_{i_k}$ be the crossings on $\alpha_i$ listed from the base point along its orientation, and let $H_c(\alpha_i) = \bigotimes\limits_{n=1}^{k} H(p_{i_n})$. Define $H_c(\beta_j) $ in a similar way. It follows that
\begin{align}
H_c = \bigotimes\limits_{i = 1}^g H_c(\alpha_i) = \bigotimes\limits_{j = 1}^g H_c(\beta_j),
\end{align}
up to a permutation of tensor components. Define
\begin{align}
\Delta_i:  H(\alpha_i) &\longrightarrow H_c(\alpha_i), & \quad &  x      \longmapsto x^{(1)} \otimes \cdots \otimes x^{(k)} \\
M_j : H_c(\beta_j) &\longrightarrow H(\beta_j),& \quad & x_1 \otimes \cdots \otimes x_k \longmapsto x_1 \cdots x_k \\
C_n : H(p_n) &\longrightarrow H(p_n),& \quad  & x \longmapsto S^{\theta(p_n)}T^{\phi(p_n)}(x)
\end{align} 
Then $\Kup{X,b;H}$ is defined by
\begin{align}
\Kup{X,b;H} = (\bigotimes\limits_{j=1}^g \mu_{\phi^j} \circ M_j)(\bigotimes\limits_{n=1}^m C_n)(\bigotimes\limits_{i=1}^g \Delta_i(e_{\phi_i})).
\end{align}

\subsection{$\hkr$ Invariant $Z_{\text{HKR}}$}
\label{subsec:hkr}
For a finite dimensional unimodular ribbon Hopf algebra $(H, R, v)$ with certain non-degeneracy condition, a topological invariant of closed oriented $3$-manifolds was constructed by Hennings \cite{Hennings} and later reformulated by Kauffman and Radford \cite{Kauffman-Radford}.
 
Given a non-zero right integral $\mu_R \in H^*$, one can associate a regular isotopy invariant $\langle L \rangle_{H,\mu_R}$ to a framed unoriented link $L$ as follows. Choose a link diagram of $L$ (still denoted by $L$) with respect to a height function such that the crossings are not critical points. On each component $L_i$ of $L$, pick a base point which is neither a crossing nor an extremum, and arbitrarily orient $L_i$. Define $\delta_i$ to be $0$ if the orientation of $L_i$ near the base point is downwards and $1$ otherwise. For a point $p$ on $L_i$ which is not an extremum, let $w_p$ be the algebraic sum of extrema between the base point and $p$, where an extremum is counted as $+1$ (resp. $-1$) if the orientation near it is counterclockwise (resp. clockwise). Equivalently, $w_p$ is $2$ times the total counterclockwise rotation, in unit of $1 = 360^{\circ}$, of the tangent of $L_i$ from the base point to $p$. Define $w_i$ to be $\frac{w_{p}}{2}$ for $p$ very close to the base point in the backward direction of $L_i$. Clearly, $w_i$ is equal to the winding number of $L_i$. Decorate each crossing with the tensor factors of the
$R$-matrix $R=\sum\limits_{i} s_i\otimes t_i$ as below.\footnote{After a crossing is decorated by the $R$-matrix elements, the over/under crossing information becomes irrelevant and we sometimes simply replace it by a solid crossing. But this is only a notation preference. In the tensor network formulation below, we will still keep the crossing as it is.}
 
 \vspace{.1in}
 
 \begin{tikzpicture}\centering
 \draw (1,1)--(-1,-1);
 \draw [color=white,line width=2mm] (-1,1)--(1,-1);
 \draw (-1,1)--(1,-1);
 \node (=) at (2,0.2) {$\leftrightarrow~\sum\limits_{i}$};
 \begin{scope}[xshift = 4cm]
 \draw (-1,1)--(1,-1);
 \draw (1,1)--(-1,-1);
 \fill[black, opacity=1] (-0.5,0.5) circle (1.5pt) node[anchor=east] {$s_i$};
 \fill[black, opacity=1] (0.5,0.5) circle (1.5pt) node[anchor=west] {$t_i$};
 \end{scope}
 \node (=) at (6,-0.7) {$,$};
 \begin{scope}[xshift = 8cm]
 \draw (-1,1)--(1,-1);
 \draw [color=white,line width=2mm] (-1,-1)--(1,1);
 \draw (-1,-1)--(1,1);
 \node (=) at (2,0.2) {$\leftrightarrow~\sum\limits_{i}$};
 \begin{scope}[xshift = 4.4cm]
 \draw (-1,1)--(1,-1);
 \draw (1,1)--(-1,-1);
 \fill[black, opacity=1] (-0.5,-0.5) circle (1.5pt) node[anchor=east] {$S(s_i)$};
 \fill[black, opacity=1] (0.5,-0.5) circle (1.5pt) node[anchor=west] {$t_i$};
 \end{scope}
 \end{scope}
 \end{tikzpicture}
 
 \vspace{.1in}
\noindent  Then we replace each decorating element $x$ on $L_i$ by $S^{-w_{p(x)}+ \delta_i}(x)$, where $p(x)$ denotes the point on $L_i$ where $x$ is located. See below for the contribution of each extremum to the powers of $S$. 
 \[
 \begin{tikzpicture}[scale=0.8]
 \begin{scope}[xshift=-6cm]
 \draw  (1,0) arc(0:180:1);
 \draw [thick,->] (0.1, 1)--(0,1);
 \node (S) at (0,1.4) {$S^{-1}$};
 \end{scope}
 \begin{scope}[xshift=-2cm]
 \draw  (1,1) arc(0:-180:1);
 \draw [thick,->] (-0.1, 0)--(0,0);
 \node (S) at (0,-0.4) {$S^{-1}$};
 \end{scope}
 \begin{scope}[xshift=2cm]
 \draw  (1,0) arc(0:180:1);
 \draw [thick,->] (0, 1)--(0.1,1);
 \node (S) at (0,1.4) {$S$};
 \end{scope}
 \begin{scope}[xshift=6cm]
 \draw  (1,1) arc(0:-180:1);
 \draw [thick,->] (0, 0)--(-0.1,0);
 \node (S) at (0,-0.4) {$S$};
 \end{scope}
 \end{tikzpicture}
 \]
 
\noindent Then $\langle L \rangle_{H,\mu_R}$ is the evaluation of the right integral $\mu_R$ on the products along each $L_i$:
\begin{equation}
\langle L \rangle_{H,\mu_R}:=\sum\limits_{(R)}\mu_R(q_1G^{1-w_1})\cdots \mu_R(q_{c(L)}G^{1-w_{c(L)}}),
\end{equation}
 where $c(L)$ is the number of
 components of $L$, $q_i\in H$ is the product of the decorating elements (after applying $S$-powers) on $L_i$ multiplied in the order following its orientation starting from the base point. 
 
 
It can be checked that $\langle L \rangle_{H,\mu_R}$ is independent of the choice of base points, orientation, and the height function. It is also preserved under framed Reidemeister moves. Thus $\langle\; \cdot\; \rangle_{H,\mu_R}$ defines an invariant of framed links.
 \begin{remark}
 \begin{enumerate}
  \item The notation here is different from but essentially the same as the Kauffman and Radford's version where the decorating elements are pushed to a vertical portion and multiplied together from bottom to top.
  \item Since $(S \otimes S)(R) = R$ and $\mu_R \circ S^2 = \mu$, one can also replace $\delta_i$ with $1-\delta_i$ in the definition of $\langle L \rangle_{H,\mu_R}$. However, this replacement has to be performed on all components of $L$ simultaneously. 
  \item If we restrict to the class of even framed links, namely, framed links where each component has an even framing, it can be shown that in any diagram of such links the winding number of each component is odd. Noting that $G^2 = u S(u^{-1})$, $\langle L\rangle_{H,\mu_R}$ can be rewritten as 
  $$\langle L \rangle_{H,\mu_R}:=\sum\limits_{(R)}\mu_R\left(q_1(u S(u^{-1}))^{\frac{1-w_1}{2}}\right)\cdots \mu_R\left(q_{c(L)}(u S(u^{-1}))^{\frac{1-w_{c(L)}}{2}}\right).$$
Hence, $\langle L\rangle_{H,\mu_R}$ does not depend on the ribbon structure of $H$ and can be defined for any unimodular $\qt$ Hopf algebras. See \cite{sawin2002invariants}.
  
 \end{enumerate}
 \end{remark}
 
\begin{figure}
\centering
\includegraphics[scale=1]{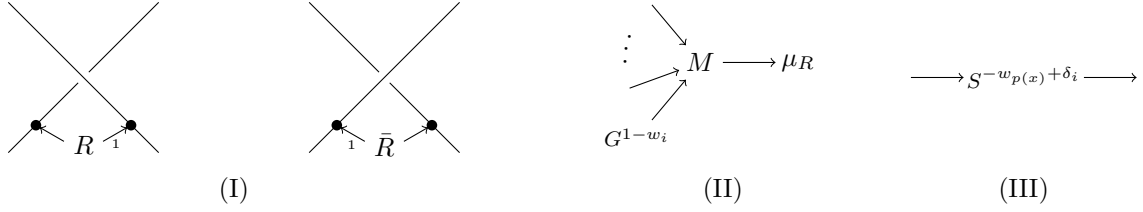}
\caption{Tensors associated with a link diagram, where $\bar{R}$ means $R^{-1}$.}\label{fig:hennings}
\end{figure}

Equivalently, it is convenient to describe $\langle L \rangle_{H, \mu_R}$ in the language of tensor networks. Again choose a base point and an orientation for each component. To each crossing assign an $R$-tensor according to the rule in Figure \ref{fig:hennings} ($\text{I}$). The first leg of the $R$- or $R^{-1}$-matrix always corresponds to the over-crossing strand. The legs terminate at links with a dot (see Figure \ref{fig:hennings} ($\text{I}$). To each component $L_i$ assign an $\widetilde{M}$-tensor as shown in Figure \ref{fig:hennings} ($\text{II}$), one leg for each dot on $L_i$ listed from the base point along its orientation. At each dot of $L_i$, insert an $S$-tensor as shown in Figure \ref{fig:hennings} ($\text{III}$) connecting the leg from the $R$-tensor to the leg from the $\widetilde{M}$-tensor. Then $\langle L \rangle_{H,\mu_R}$ is equal to the contraction of these tensors.  

It is a direct calculation that the invariant of the unknot with framing $\pm 1$ is $\mu_R(v^{\pm 1})$. From now on assume $\mu_R(v)\mu_R(v^{-1})\neq 0$, which is the non-degeneracy condition we impose on $H$ and which is always true when $H$ is factorizable \cite{Cohen-Westreich}. Let $\omega(v)$ be a square root of $\mu_R(v)/\mu_R(v^{-1})$, then $\mu_R(v)/\omega(v)$ is a square root of $\mu_R(v) \mu_R(v^{-1})$. The $Z_{\text{HKR}}$ invariant for a closed oriented $3$-manifold $X$ is defined to be:
 \begin{equation}
 \HKR{X;H, \omega(v)}=(\mu_R(v)/\omega(v))^{-c(L)}\;\omega(v)^{-\sign(L)}\;\langle L \rangle_{H, \mu_R},
 \end{equation}
 where $L$ is a surgery link of $X$ and $\sign(L)$ denotes the signature of the framing matrix of $L$.

 \begin{remark}

By its very definition, $\omega(v)$ does not depend on $\mu_R$. For any non-zero scalar  $s \in \Complex$, clearly we have $\langle L \rangle_{H, s \mu_R} = s^{c(L)} \langle L \rangle_{H,\mu_R}$. It follows that $\HKR{X;H, \omega(v)}$ does not depend on $\mu_R$ either. If one chooses the other square root $-\omega(v)$, then 
\begin{align*}
\HKR{X;H, \omega(v)} = (-1)^{c(L) + \sign(L)} \HKR{X;H, -\omega(v)}.
\end{align*}
Hence, up to a negative sign $\HKR{X;H, \omega(v)}$ does not depend on the choice of a square root of $\mu_R(v)/\mu_R(v^{-1})$, in which case the invariant is more commonly written as:
\begin{equation*}
 \HKR{X;H}=[\mu_R(v)\mu_R(v^{-1})]^{-\frac{c(L)}{2}}\;[\mu_R(v)/\mu_R(v^{-1})]^{-\frac{\sign(L)}{2}}\;\langle L \rangle_{H,\mu_R}.
 \end{equation*}

 \end{remark}
 
Just as the $\wrt$ invariant, $Z_{\text{HKR}}$ can also be refined to an invariant of $3$-manifolds endowed with a $2$-framing. We recall the definition of a $2$-framing introduced in \cite{atiyah1990framings}. Let $N$ be a Riemannian manifold of dimension $n \geq 3$. Consider the diagonal embedding of $\SO(n)$ into $\SO(2n)$: 
\begin{equation*}
\centering
\begin{tikzcd}
 & & \Spin(2n) \arrow[d, " \pi"] \\
\SO(n) \arrow[r, " "] \arrow[urr,  dashed, " "] & \SO(n) \times \SO(n) \arrow[r, " "] & \SO(2n) 
\end{tikzcd}
\end{equation*}
The embedding induces a lift from $\SO(n)$ to $\Spin(2n)$, indicated by the dashed arrow, so that the diagram above commutes. The diagram determines a spin structure of $2T_N:= T_N \oplus T_N$, double of the tangent bundle of $N$. A $2$-framing of $N$ is defined to be  a  trivialization of $2T_N$ viewed as a $\Spin(2n)$ bundle. For three manifolds, $2$-framings are equivalent to $p_1$ structures \cite{blanchet1995topological}. Let $X$ be a closed oriented $3$-manifold. Since $\pi_1(\Spin(6)) = \pi_2(\Spin(6)) = 0, $ $\pi_3(\Spin(6)) = \Integer$, the set of $2$-framings of $X$ form a torsor over $H^3(X; \pi_3(\Spin(6))) \simeq \Integer$. Choose any $4$-manifold $W$ whose boundary is $X$. For any $2$-framing $\phi$ on $X$, define
\begin{align}
\sigma(\phi):= 3\sign(W) - \frac{1}{2} p_1(2T_W, \phi),
\end{align}
where $\sign(W)$ is the Hirzebruch signature of $W$ and $p_1(2T_W, \phi)$ is the relative Pontrjagin number. \footnote{Note that the $\sigma$ map defined here is equal to three times the $\sigma$ invariant in \cite{atiyah1990framings}.} By the Hirzebruch signature formula for closed $4$-manifolds, $\sigma(\phi)$ is independent of the bounding manifold $W$. Since $2T_X$ is spin, it implies $p_1(2T_W, \phi)$ is an even integer. Moreover, $\sigma$ is an affine linear isomorphism from the set of the $2$-framings to $\Integer$. The canonical $2$-framing is the unique  $\phi_0$ satisfying $\sigma(\phi_0) = 0$.

Let $H, \mu_R, v$ be as above, $\omega_6(v)$ be a sixth root of $\mu_R(v)/\mu_R(v^{-1})$ and $\omega(v) = \omega_6(v)^3$. The $Z_{\text{HKR}}$ invariant for the pair $(X, \phi)$ is defined to be:
 \begin{equation}
 \HKR{X, \phi; H, \omega_6(v)}:=(\mu_R(v)/\omega(v))^{-c(L)}\;\omega_6(v)^{-\frac{1}{2}p_1(2T_{W_L}, \phi)}\;\langle L \rangle_{H, \mu_R},
 \end{equation}
 where $W_L$ is the $4$-manifold obtained from the surgery link $L$. It follows immediately from the definitions that
 \begin{align}
 \HKR{X, \phi; H, \omega_6(v)} = \omega_6(v)^{\sigma(\phi)} \HKR{X; H, \omega(v)}.
 \end{align}
Thus the original invariant is equal to the refined invariant evaluating at the canonical $2$-framing. The chosen roots $\omega_6(v) $ and $\omega(v)$ are often dropped from the formula when they are clear from the context. In the following we use $\HKR{\cdot}$ to denote both the refined invariant and the original one.

\section{Main Results I}
\label{sec:main1}
In this section, $H$ denotes a finite dimensional $\db$ Hopf algebra. Hence its $\Drinfeld$ double $D(H) = H^{*\cop} \otimes H$ is ribbon. Note that $D(H)$ is also factorizable and unimodular. See Section \ref{sec:Hopf} for our notations on Hopf algebras.
\begin{theorem}[$=$ Theorem \ref{thm:main11} ]
\label{thm:main1}
Let $H$ be a finite dimensional $\db$ Hopf algebra and $X$ be a closed oriented $3$-manifold, then there exist a framing $b$ and a $2$-framing $\phi$ of $X$ such that,
\begin{align}
\Kup{X,b;H} = \HKR{X,\phi; D(H)}.
\end{align}
\begin{proof}
The proof is given in the next three subsections. Section \ref{subsec:heegaard surgery} gives a special Heegaard diagram of $X$ in which one family of circles form a surgery link for $X$. In Section  \ref{subsec:computing Kuperberg}, we construct a framing $b$ of $X$ presented on the Heegaard diagram and compute $\Kup{X,b;H}$. In Section \ref{subsec:computing hkr} we define a $2$-framing $\phi$ and compute $\HKR{X, \phi; D(H)}$. The equality in the theorem then follows.
\end{proof} 
\end{theorem}
 \subsection{Special Heegaard Diagrams}
\label{subsec:heegaard surgery}
A Heegaard diagram is a triple $R = (\Sigma_g, \alpha, \beta)$ where $\Sigma_g$ is a closed oriented surface of genus $g$, and $\alpha = \{\alpha_1, \cdots, \alpha_g\}$ (resp. $\beta = \{\beta_1, \cdots, \beta_g\}$) is a collection of $g$ disjoint simple closed curves such that the complement of the $\alpha_i\;'$s (resp. the $\beta_j\;'$s) in $\Sigma_g$ is a $2g$-punctured sphere. A closed oriented $3$-manifold is obtained from a Heegaard diagram by attaching $2$-handles to the closed curves and filling sphere boundaries with $3$-handles. Every closed oriented $3$-manifold can be represented by a Heegaard diagram, and different diagrams of the same manifold are related by isotopy, handle slides and stabilization. A diagram $R = (\Sigma_g, \alpha, \beta)$ of the $3$-sphere $\Sphere^3$ is {\it standard} if the geometric intersection of $\alpha_i$ with $\beta_j$ is $1$ for $i=j$ and $0$ otherwise. 
Every standard diagram of genus $g$ for $\Sphere^3$ is isotopic to the one obtained by taking stabilization $g$ times from the two sphere.
Heegaard diagrams with certain special properties are studied in \cite{birman1977special}.

\begin{theorem}\cite{birman1977special}
\label{thm:special}
Every closed oriented $3$-manifold $X$ has a Heegaard diagram $R = (\Sigma_g, \alpha, \beta)$ for some genus $g$ satisfying the following properties:
\begin{enumerate}
\item There exists a collection of $g$ curves $\gamma = \{\gamma_1, \cdots, \gamma_g\}$ on $\Sigma_g$ such that both $R_1 = (\Sigma_g, \alpha, \gamma)$ and $R_2 = (\Sigma_g, \beta,\gamma)$ are standard diagrams for $\Sphere^3$.
\item View $\beta$ as a framed link in $\Sphere^3$ determined by $R_1$, where the framing is taken to be a parallel copy of $\beta$ in the Heegaard surface. Then $\beta$ is a surgery link for $X$. Moreover, the framings are all even integers. 
\end{enumerate}
\begin{proof}
(Sketch) See \cite{birman1977special} for a more detailed proof. It is a standard result that $X$ has a surgery link $L$ which is the plat closure of a certain $2g$-strand braid $\sigma \in B_{2g}$. Actually one can always choose $\sigma$ to be a pure braid and the framing of each component to be an even integer. In this case,  $L$ has $g$ components $\{L_1, \cdots, L_g\}$. 
Assume $\sigma$ is aligned vertically in the stripe $\{0\} \times \Real \times [0,1]$ with end points $(0, i, 0), \, (0, i, 1)$, $i = 1,2, \cdots, 2g$. The $i$-th plat on the bottom (resp. on the top) connects $(0,2i-1,0)$ and $(0,2i,0)$ (resp. $(0,2i-1,1)$ and $(0,2i,1)$). See Figure \ref{fig:braid} (Left). According to Theorem $5.2$ in \cite{birman1977special}, one can isotope $L$, by untwisting the braid at the cost of twisting the plats on the top, so that each $L_i$ is decomposed as $ L_i^1 \cup L_i^2$ (see Figure \ref{fig:braid} (right)), where $L_i^1$ is the arc consisting of the segments $\{0\} \times \{2i-1, 2i\} \times [0,1]$ and the $i$-th plat on the bottom, and $L_i^2$ is an arc in $\Real^2 \times \{1\}$ connecting $(0,2i-1,1)$ and $(0,2i,1)$. Moreover, the $L_i^2\;'$s are disjoint from each other. Arbitrarily choose $g-1$ mutually disjoint arcs $C_1, \cdots, C_{g-1}$ in the plane $\Real^2 \times \{1\}$ so that $C_j$ connects a point in $L_j^2$ to a point in $L_{j+1}^2$ and is otherwise disjoint from all the $L_i^2\;'$s. Let
\begin{align*}
B:= \left(\bigcup\limits_{i=1}^{g-1} C_i\right) \bigcup \left(\bigcup\limits_{i=1}^{g}L_i^2\right), \qquad H:= B \bigcup \left(\bigcup\limits_{i=1}^{g} L_i^1\right),
\end{align*}
and $N(B)$ and $N(H)$  be a regular neighborhood of $B$ and $H$, respectively. Then $N(B)$ is a $3$-ball and $N(H)$ is a handlebody obtained from $N(B)$ by attaching $g$ $1$-handles, each of which corresponds to a regular neighborhood $N(L_i^1)$ of $L_i^1$. Clearly $\Sphere^3 = N(H) \cup N(H)^c$ is a Heegaard decomposition of $\Sphere^3$. On $\partial N(H)$ choose a complete set of meridian curves $\gamma = \{\gamma_1, \cdots, \gamma_g\}$ for $N(H)$ and a complete set of meridian curves $\alpha = \{\alpha_1, \cdots, \alpha_g\}$ for $N(H)^c$ so that $(\partial N(H), \alpha, \gamma)$ is a standard Heegaard diagram of $\Sphere^3$. 

Let $\beta = \{\beta_1, \cdots, \beta_g\}$ be a set of curves with $\beta_i$ representing the framing of $L_i$. 
One can assume $\beta_i$ is contained in $\partial N(H) \cap \partial N(L_i)$, where $N(L_i)$ is a regular neighborhood of $L_i$. It follows that the complement of $\beta$ in $\partial N(H)$ is a $2g$-punctured sphere. It can be shown that $(\partial N(H), \beta, \gamma)$ is a standard Heegaard diagram of $\Sphere^3$ and that $(\partial N(H), \alpha, \beta)$ is a Heegaard diagram of $X$. Clearly $L$ and $\beta$ are isotopic framed links.
\begin{figure}
\centering
\includegraphics[scale=1]{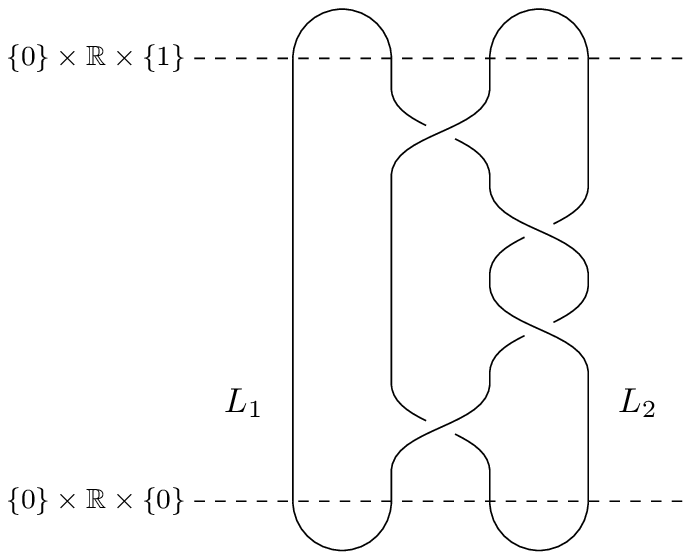}
\qquad\quad
\includegraphics[scale=1]{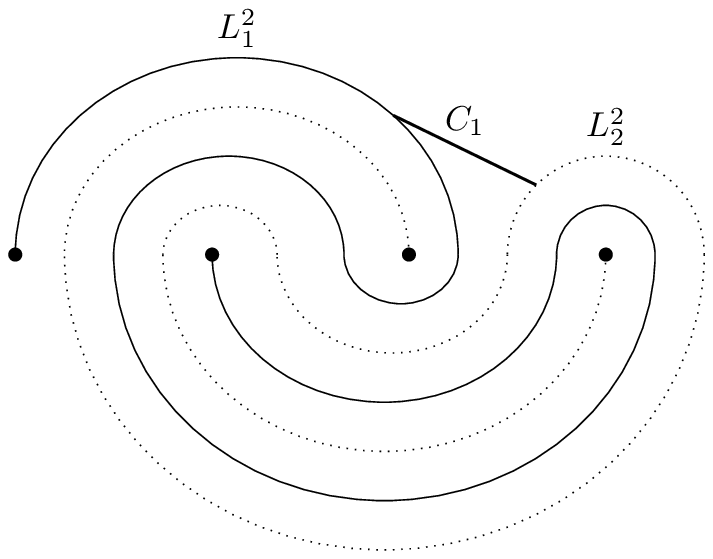}
\caption{(Left) $L$ is the plat closure of $\sigma = \sigma_2 \sigma_3^{-2} \sigma_2$; (Right) The part of $L$ in the plane $\Real^2 \times \{1\}$ where $L_i^1$ (not drawn) is an arc lying inside the page connecting the $(2i-1)$-th dot and the $(2i)$-th dot, $i = 1,2$.   }\label{fig:braid}
\end{figure}
\end{proof}
\end{theorem}  


Theorem \ref{thm:special} provides a bridge between Heegaard diagrams and surgery links which is exactly the ingredient that will be used to compare $Z_{\text{Kup}}$ and $Z_{\text{HKR}}$. For the sake of clarity, we give an explicit description of the Heegaard diagram/surgery link model.

\begin{figure}
\centering
\includegraphics[scale=1]{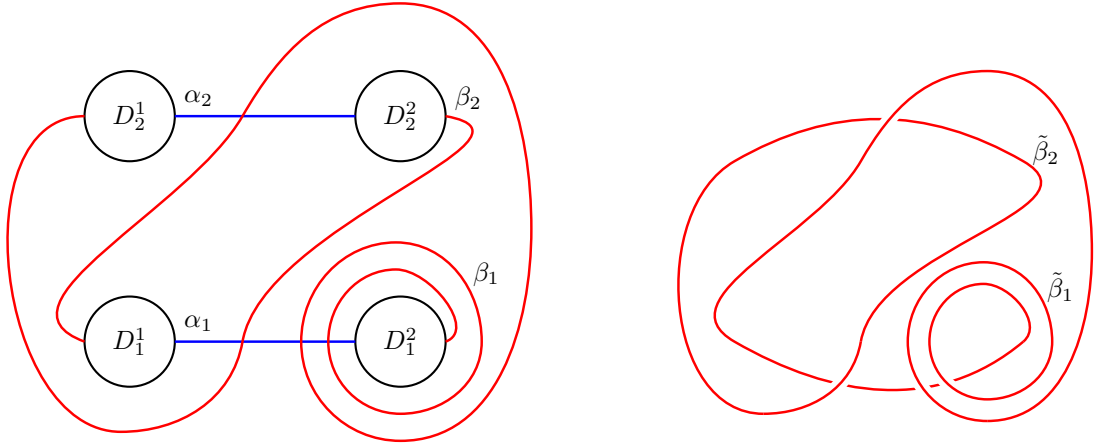}
\caption{(Left) The $x$-$y$ plane. $\alpha$ and $\beta$ are represented by arcs in blue and red, respectively. The $1$-handles $N_i\;'$s are outside of the plane. (Right) The diagram $\tilde{\beta}$ in $x$-$y$ plane of $\beta$ viewed as a link.} \label{fig:special_diagram}
\end{figure}

Endow $\Real^3$ with the $\{x,y,z\}$ coordinates. Let $\Sphere^3 = \Real^3 \cup \{\infty\}$, $B_{+} = ( \Real^2 \times \Real_{\geq 0}) \cup \{\infty\}$, $B_{-} = ( \Real^2 \times \Real_{\leq 0}) \cup \{\infty\}$, and $\Sphere^2 = B_{+} \cap B_{-} =  (\Real^2 \times \{0\}) \cup \{\infty\}$. Also identify $\Real^2 $ with $\Real^2 \times \{0\} \subset \Sphere^2$. Fix an integer $g \geq 0$. For $1 \leq i \leq g$, let $D_i^1$ and $D_i^2$ be the disks in $\Sphere^2$ centered at $(0,i)$ and $(1,i)$, respectively, of radius $\epsilon \ll 1/8$, and let $N_i$ be a three dimensional $1$-handle in $B_{-}$ connecting $D_i^1$ and $D_i^2$. The $N_i\;'$s are unknotted and unlinked.
For instance, one can push the segment $ [0,1] \times \{i\}$ slightly into $B_{-}$ keeping the end points fixed and set $N_i$ to be a regular neighborhood of the push-off.
Set $B_{g,+} = B_{+} \cup_{i=1}^{g} N_i$ and $B_{g,-} = B_{-} \setminus (\cup_{i=1}^{g} N_i)$, then $\Sphere^3 = B_{g,+} \cup B_{g,-}$ is a standard genus-$g$ Heegaard decomposition. Define $\Sigma_g := \partial B_{g,+}= \partial B_{g,-}$, $\partial N_i:= \Sigma_g \cap N_i$, and $\Sphere^{2;2g}:=   \Sigma_g \cap \Sphere^2$. Clearly $\Sphere^{2;2g}$ is a $2g$-punctured sphere. We call $\partial D_i^1$ and $\partial D_i^2$ the left foot and right foot, respectively, of $\partial N_i$. The readers may find Figure \ref{fig:special_diagram} helpful in the following discussions. Take $\alpha_i$ to be a meridian of $B_{g,-}$ which consists of the segment $ [\epsilon, 1-\epsilon] \times \{i\}$ and the arc traveling through $\partial N_i$ once (without twisting around $N_i$) connecting $(\epsilon,i)$ and $(1-\epsilon,i)$. Also take $\gamma_i$ to be a meridian of $B_{g,+}$ circling a section of $N_i$ once. Let $\beta_i$ be any simple closed curve in $\Sigma_g$ which travels through $\partial N_i$ once (without twisting around $N_i$) and spends the rest of time in $\Sphere^{2;2g}$. Moreover, $\beta_i$ is parallel to $\alpha_i$ when traveling in $\partial N_i$ and all the $\beta_i\;'$s are disjoint from each other. Furthermore, each $\beta_i$ crosses the segment $ [\epsilon, 1-\epsilon] \times \{i\}$ an even number of times.  Define $\alpha = \{\alpha_1,\cdots, \alpha_g\}$ and define $\beta, \gamma$ analogously. Since $\beta$ is contained in $\Sigma_g$, it can be naturally viewed as a framed link in $\Sphere^3$ by taking a parallel copy in $\Sigma_g$ as the framing curve. Furthermore, it has a diagram in $\Sphere^2$ obtained by projecting the part of each $\beta_i$ in $\partial N_i$ to $[0,1] \times \{i\}$ while keeping the part in $\Sphere^{2;2g}$ fixed. See Figure \ref{fig:special_diagram} (Right). Denote the projection by $\tilde{\beta}$. With notations from above and by Theorem \ref{thm:special}, we have
\begin{enumerate}
  \item $(\Sigma_g,\alpha,\gamma)$ is a standard Heegaard diagram of $\Sphere^3$.
  \item $\tilde{\beta}$ is a link diagram for $\beta$, and the self-linking number of each component of $\tilde{\beta}$ is an even integer.
  \item $(\Sigma_g,\alpha,\beta)$ is a Heegaard diagram for the $3$-manifold whose surgery link $\beta$. Denote such a $3$-manifold by $X(\Sigma_g,\beta)$ with $\alpha$ and $\gamma$ known implicitly. Then every closed oriented $3$-manifold is homeomorphic to some $X(\Sigma_g,\beta)$ .
\end{enumerate}

\subsection{A Framing on $X(\Sigma_g,\beta)$ and the Kuperberg Invariant}
\label{subsec:computing Kuperberg}
Given the $3$-manifold $X = X(\Sigma_g, \beta)$, we construct a framing of $X$ presented in the Heegaard diagram $(\Sigma_g, \alpha, \beta)$. Recall from Section \ref{subsec:Kuperberg} that a framing consists of two orthogonal combings $b_1$ and $b_2$ satisfying certain conditions, where $b_1$ is represented as a vector field with $2g+1$ singularities and $b_2$ is represented as a set of twist fronts. For $1 \leq i \leq g$, let $w_i$ be the winding number of $\tilde{\beta}_i$. Since the framing of $\tilde{\beta}_i$ is even, then $w_i$ is odd. Set $w_i = 2n_i + 1$. 

\begin{figure}
\centering
\includegraphics[scale=1]{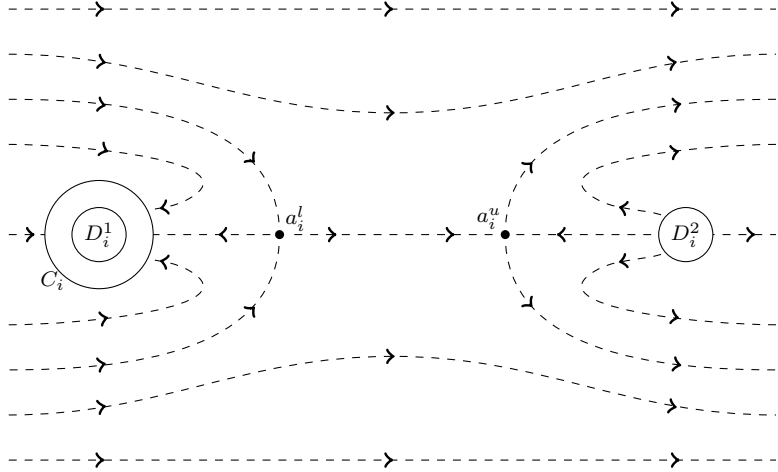}
\caption{Flow lines of $b_1$ in the closed rectangle $\overline{R_i} = [-1,2] \times [i-1/2, i+1/2] \subset \Real^2$.}\label{fig:framing_b1}
\end{figure}

\begin{figure}
	\centering
	\includegraphics[scale=1]{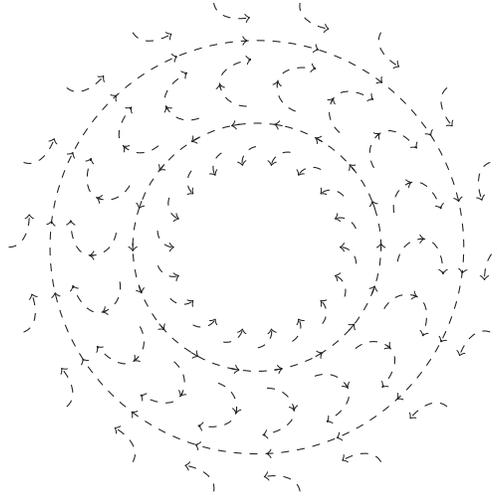}
	\caption{A model of $b_1$ rotation in the annulus between $C_i$ and $\partial D_i^1$ for $n_i = 1$. The two circles are closed flow lines, but not $C_i$ and $\partial D_i^1$.}
	\label{fig:rotation}
\end{figure}

\vspace{0.5cm}
\noindent \textbf{First combing} $\mathbf{b_1}$:  the construction of $b_1$ is generalized from that given in \cite{Chang-Wang}.  We describe the flow lines and singularities of $b_1$. The singularities are located at $a_i^l := (1/4, i),\, a_i^u = (3/4, i)$, and $\infty$, $i=1, \cdots, g$. All the singularities have index $-1$ except the one at $\infty$ which has index $2$.  Let $R_i$ be the open rectangle $(-1,2) \times (i-1/2, i+1/2)$ and $R = \sqcup_{i=1}^{g} R_i$. In $\Real^2 \setminus R$, $b_1$ takes the value $\frac{\partial}{\partial x}$, i.e., $b_1$ points toward the positive direction of the $x$-axis \footnote{It is direct to check this implies $\infty$ is a singular point of index $2$.}. Note that on the boundary of each $R_i$, the value of $b_1$ is $\frac{\partial}{\partial x}$. Now it suffices to describe $b_1$ inside $R_i$ and $\partial N_i$. This is illustrated in Figure \ref{fig:framing_b1}, where dashed lines represent the flow lines and $C_i$ is the circle centered at $(0,i)$ with radius $2\epsilon$. The behavior of $b_1$ inside the annulus bounded by $C_i$ and $\partial D_i^1$ is as follows. The field $b_1$ points toward the center on $C_i$ and $\partial D_i^1$. Along each radial segment connecting $C_i$ and $\partial D_i^1$, $b_1$ rotates counterclockwise, in unit $1 = 360^{\circ}$, by the degree $n_i$.\footnote{If $n_i < 0$, then the rotation is clockwise of degree $-n_i$.} If we set the center of $C_i$ to be $(0,0)$ for simplicity, then a formula of $b_1$ inside the annulus is given by:
\begin{align}
b_1(x,y) = -\cos\left(\theta + 2\pi n_i \frac{2\epsilon-r}{\epsilon}\right)\frac{\partial}{\partial x} - \sin\left(\theta + 2\pi n_i \frac{2\epsilon-r}{\epsilon}\right)\frac{\partial}{\partial y},
\end{align}
where $(r,\theta)$ is the polar coordinate of $(x,y)$. Note that the radial segments are not flow lines.  Figure \ref{fig:rotation} shows a model of flow lines for $n_i = 1$. The rotation of any degree can be obtained by stacking this model or the orientation reversal model in the radial direction. Inside the tube $\partial N_i$ the flow lines of $b_1$ travel from one end to the other without any twisting and emerges out of $\partial D_i^2$.

\begin{figure}
\centering
\includegraphics[scale=1]{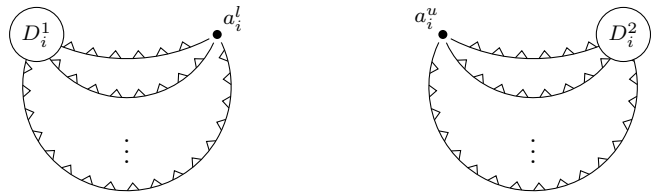}
\caption{Twist fronts of $b_2$ connecting $a_i^l$ and $a_i^u$ in the case $w_i > 0$.}\label{fig:framing_b2}
\end{figure}

\vspace{0.5cm}
\noindent \textbf{Second combing} $\mathbf{b_2}$: for $1 \leq i \leq g$, there are $|w_i|$ twist fronts, each of which travels through $\partial N_i$ in parallel and connects the two singularities $a_i^l$ and $a_i^u$. See Figure \ref{fig:framing_b2}. The (small triangles on) twist fronts point upward as shown in the figure if $w_i > 0$ and downward otherwise. 

\begin{figure}
\centering
\includegraphics[scale=1]{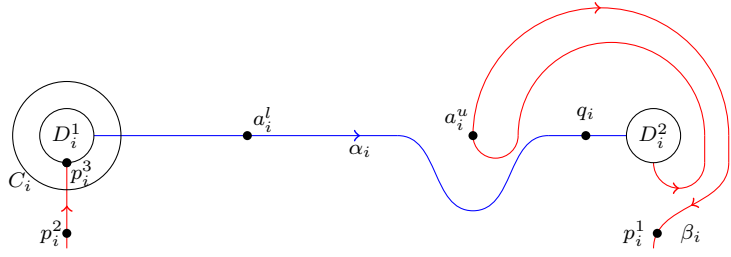}
\caption{Perturbations, base points, and orientations of lower and upper circles}\label{fig:lower_upper_perturb}
\end{figure}

\vspace{0.5cm}
\noindent \textbf{Lower and upper circles}: we designate $\alpha$ and $\beta$ as the set of lower and upper circles, respectively. But note that we need each circle to pass exactly one singular point of index $-1$ in a specific manner (see Section \ref{subsec:Kuperberg}). We achieve this by perform a slight perturbation on the circles. See Figure \ref{fig:lower_upper_perturb}. For each $i$, set the base point of $\alpha_i$ to be $a_i^l$ and orient $\alpha_i$ so that it points to the positive $x$-direction (horizontally to the right in the figure) at $a_i^l$. Then perturb $\alpha_i$ off $a_i^u$ and perturb $\beta_i$ so that it passes $a_i^u$. Set $a_i^u$ as the base point of $\beta_i$. The orientation is chosen so that it points upward at $a_i^u$.

By isotopy, we may assume that each $\beta_i$ is away from the feet of all $\partial N_j\,'$s and from all base points except the part as shown in Figure \ref{fig:lower_upper_perturb}. In particular, all the intersections of the lower circles with upper circles are constrained in the horizontal segments connecting a lower base point to the corresponding upper base point. At the intersections, the upper circles are vertical. Recall from Section \ref{subsec:Kuperberg} the definitions of $\theta_c,\, \theta_c(p), \,\phi_c, \,\phi_c(p)$. 
Let $p,q$ be two points on a lower or upper circle $c$ and define $\theta_{c}(p,q) = \theta_{c} (q) - \theta_{c}(p)$, namely, $\theta_{c}(p,q)$ is the degree of rotation of $c'$ relative to $b_1$ from $p$ to $q$ along $c$.

\begin{lemma}\label{lem:theta}
\begin{itemize}
\item Let $p_i^1, p_i^2, p_i^3$ be points on $\beta_i$ as shown in Figure \ref{fig:lower_upper_perturb}, then
\begin{align}
\theta_{\beta_i}(a_i^u, p_i^1) = \theta_{\beta_i}(p_i^1) = -\frac{1}{4}, \quad \theta_{\beta_i}(p_i^2, p_i^3) = -\frac{1}{4} - n_i, \quad \theta_{\beta_i}(p_i^3, a_i^u) = -\frac{1}{2}.
\end{align}
\item Let $q_i$ be a point on $\alpha_i$ as shown Figure \ref{fig:lower_upper_perturb}, then $\theta_{\alpha_i}(q_i) = \frac{1}{2}$.
\item Let $p$ be a point on $\beta_i$ between $p_i^1$ and $p_i^2$ and assume the tangent of $\beta_i$ at $p$ is vertical, then $\theta_{\beta_i}(p_i^1, p) = \frac{w_p}{2}$, where $w_p$ (also see Section \ref{subsec:hkr}) is the algebraic sum of extrema along $\beta_i$ between $p_i^1$ to $p$, where an extremum is counted as $+1$ if the orientation near it is counterclockwise, and as $-1$ otherwise.
\end{itemize}
\begin{proof}
The first two parts follow directly from observations of Figure \ref{fig:framing_b1} and \ref{fig:lower_upper_perturb}. In particular, $\theta_{\beta_i}(p_i^2, p_i^3)$ would be $-\frac{1}{4}$ if the flow lines inside the annulus between $C_i$ and $\partial D_i^1$ did not rotate. The rotations in the annulus by the degree $n_i$ contributes an extra $-n_i$ to $\theta_{\beta_i}(p_i^2, p_i^3)$. The third part is obtained by noting that when traveling along $\beta_i$ away from all base points, each pass of an extremum contributes $\pm \frac{1}{2}$ to $\theta_{\beta_i}(p_i^1, p)$ depending on the orientation near the extremum. Also see Lemma $1$ in \cite{Chang2015}. 
\end{proof}
\end{lemma} 

\begin{lemma}\label{lem:b1b2}
For the combings $b_1, b_2$ constructed above, we have 
\begin{align}
\theta_{\alpha_i} = \phi_{\alpha_i} = n_i + \frac{1}{2}, \quad \theta_{\beta_i} = -\phi_{\beta_i} = n_i + \frac{1}{2}, \quad 1 \leq i \leq g.
\end{align} 
\begin{proof}
By the third part of Lemma \ref{lem:theta}, we have $\theta_{\beta_i}(p_i^1, p_i^2) = w_i + \frac{1}{2}$, where $w_i$ is the winding number of $\tilde{beta}_i$. Hence, 
\begin{align*}
\theta_{\beta_i} &= \theta_{\beta_i}(a_i^u, p_i^1) + \theta_{\beta_i}(p^1, p_i^2) + \theta_{\beta_i}(p^2, p_i^3) + \theta_{\beta_i}(p^3, a_i^u) \\
                 &= (-\frac{1}{4}) + (w_i + \frac{1}{2}) + (-\frac{1}{4} -n_i) + (-\frac{1}{2}) = n_i + \frac{1}{2}.
\end{align*}

For $\theta_{\alpha_i}$, note that when traveling along $\alpha_i$ from $q_i$ to $a_i^l$, we will cross the annulus between $\partial D_i^1$ and $C_i$, and the direction of the crossing is from $\partial D_i^1$ to $C_i$. During this crossing, the vector field $b_1$ rotates by a degree of $-n_i$, and hence $\theta_{\alpha_i}$ increases by $n_i$. Then the equality $\theta_{\alpha_i} =  n_i + \frac{1}{2}$ follows from the second part of Lemma \ref{lem:theta}.

The equalities concerning the $\phi\,'s$ are derived by counting the number of crossings of the circles with twist fronts.
\end{proof}
\end{lemma}

By Lemma \ref{lem:b1b2}, the combings $b_1, b_2$ extend to a framing on $X$. Denote this framing by $b = (b_1, b_2)$. 

\begin{lemma}
Let $p$ be a crossing of $\beta_i$ with $\alpha_j$, then $\theta_{\alpha_j}(p)  = \phi_{\alpha_j}(p) = \phi_{\beta_i}(p)= 0$, and  $\theta_{\beta_i}(p) = -\frac{1}{4} + \frac{w_p}{2}$, where $w_p$ is defined as in the third part of Lemma \ref{lem:theta}. In particular, in the tensor network computing $\Kup{X,b;H}$, the tensor assigned to $p$ is $S^{\theta(p)}T^{\phi(p)}$ where $\theta(p) = 1- w_p$ and $\phi(p) = 0$.
\end{lemma}

The Kuperberg invariant $\Kup{X,b;H}$ can be described as follows. Assign the tensors in Figure \ref{fig:three_tensors} to each $\alpha_i$, each $\beta_j$, and each crossing $p$, with $\phi_i = n_i + \frac{1}{2}$, $\phi^j = -n_i - \frac{1}{2}$, $\phi(p) = 0$, and $\theta(p) = 1-w_p$.

\subsection{Computing  $Z_{\text{HKR}}$}
\label{subsec:computing hkr}
\begin{figure}
\centering
\includegraphics[scale=1]{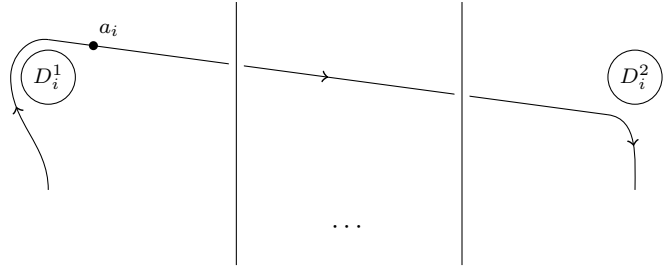}
\caption{A perturbation of the diagram $\tilde{\beta}$.}\label{fig:link_perturb}
\end{figure}

 We compute $Z_{\text{HKR}}$ for the $3$-manifold $X = X(\Sigma_g, \beta)$ from $D(H)$. See Section \ref{sec:Hopf} and  \ref{subsec:hkr} for some notations to be used below.  Recall from Section \ref{subsec:heegaard surgery} that a surgery link diagram for $X$ is $\tilde{\beta}$. We perturb $\tilde{\beta}$ slightly so that the $y$-coordinate function serves as a height function for $\tilde{\beta}$. The perturbed diagram, still denoted by $\tilde{\beta}$, is shown in Figure \ref{fig:link_perturb}. That is, instead of connecting the two feet $\partial D_i^1, \partial D_i^2$ horizontally, $\tilde{\beta}_i$ travels from slightly over the top of $D_i^1$ to slightly below the bottom of $D_i^2$ in a right-downwards direction. We also assume all the crossings are right-handed and are constrained in the segments $\sqcup_{i=1}^g (1/4,3/4) \times \{i\}$. 
Pick a point $a_i$ on $\tilde{\beta}_i$ near the left feet $\partial D_i^1$ (past the maximum ) as the base point of $\tilde{\beta}_i$ and orient $\tilde{\beta}_i$ so that it points to the right feet $\partial D_i^2$ at $a_i$. Under this orientation, we have $\delta_i = 0$.  

\begin{figure}
\centering
\includegraphics[scale=1]{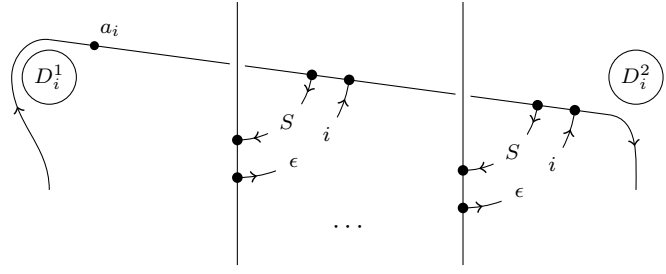}
\caption{An assignment of $R$ tensor to each crossing}\label{fig:R_assignment}
\end{figure}

\begin{figure}
\centering
\includegraphics[scale=1]{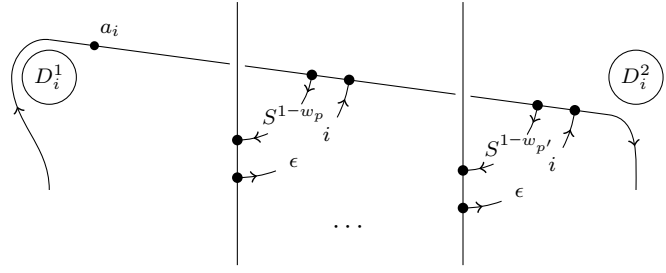}
\caption{Combining the $R$ tensor and $S$-tensor.}\label{fig:R_S_assignment}
\end{figure}

\begin{figure}
\centering
\includegraphics[scale=1]{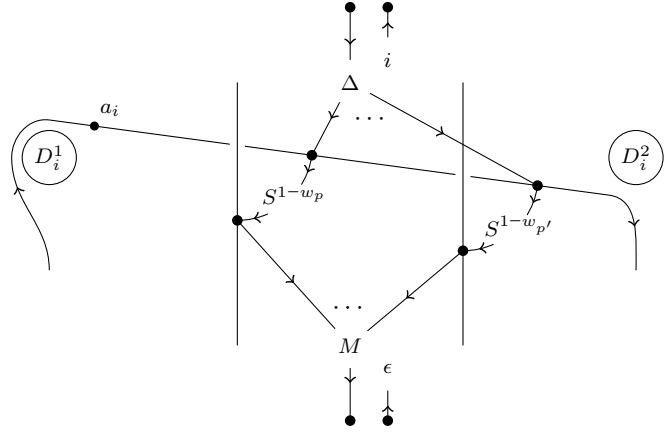}
\caption{Assignment of $\widetilde{M}$ tensor}\label{fig:M_assignment1}
\end{figure}

\begin{figure}
\centering
\includegraphics[scale=1]{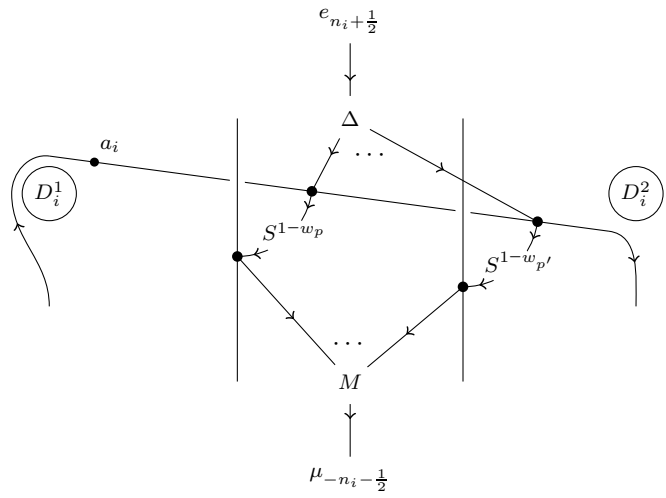}
\caption{Contraction of the $R$, $S$, and $\widetilde{M}$ tensors}\label{fig:M_assignment2}
\end{figure}

We use tensor network formulation to compute $\langle \tilde{\beta} \rangle_{D(H), \mu_R^D}$. Recall that in $D(H)$, each leg in a tensor consists of two lines, one corresponding to $H^{*\cop}$ and the other to $H$. The $(R^D)^{-1}$ tensor is assigned to all crossings since they are all right-handed. See Figure \ref{fig:R_assignment}. A dot at the end of a leg indicates a position where tensors will be contracted later. Note that here two neighboring dots are treated as one dot since we are working with tensors in $D(H)$. Call a dot {\it covariant} if the leg attached to it is incoming and {\it contravariant} otherwise. We examine the $S$-tensor assigned to each dot. Each dot on a horizontal segment has an $S^D$ power of $0$ since there are no extrema between the base point to where the dot is located. For a dot on a vertical segment corresponding to a crossing $p$, assume it belongs to some $\tilde{\beta}_i$, then its $S_D$ power is $-w_p$ where $w_p$ is the algebraic sum of extrema between $a_i$ and the dot.  Note that $S_D^{-w_p}(\epsilon \otimes x) = \epsilon \otimes S^{-w_p}(x)$. Combining the $R$-tensor and $S$-tensor, the configuration now is as in Figure \ref{fig:R_S_assignment}. Finally we apply the $\widetilde{M}$ tensor in Figure \ref{fig:hennings} to each $\tilde{\beta}_i$. This is broken down to several stages. Start from the base point $a_i$ and travel along $\tilde{\beta}_i$ following its direction. One first comes across dots on the horizontal segment, and then dots on vertical segments. Firstly, multiplying the elements on the horizontal segments is equivalent to attaching a $\Delta$-type tensor in Equation \ref{equ:general_Delta_M} (Left)  with each outgoing leg corresponding to a contravariant dot from left to right. Secondly, multiplying elements on the vertical segments is equivalent to attaching an $M$-type tensor in Equation \ref{equ:general_Delta_M} (Right) with each incoming leg corresponding to a covariant dot, which again corresponds to the crossings on $\tilde{\beta}_i$. See Figure \ref{fig:M_assignment1}. Recall that $w_i = 2n_i + 1$ is the winding number of $\tilde{\beta}_i$. Finally, the whole $\widetilde{M}$-tensor is obtained by multiplying the two dots on the top (Figure \ref{fig:M_assignment1}), the two dots on the bottom (Figure \ref{fig:M_assignment1}), and the element $ (a^D)^{-n_i} = \alpha^{n_i} \otimes a^{-n_i}$, followed by the application of $\mu_R^D = e_L \otimes \mu_R$. Note that for $f \in H^*, x \in H$,
\begin{align*}
\mu_R^D\left((f \otimes i)(\epsilon \otimes x)(\alpha^{n_i} \otimes a^{-n_i})\right) = \mu_R^D\left((f \otimes x)(\alpha^{n_i} \otimes a^{-n_i})\right) = f(e_{n_i + \frac{1}{2}})  \mu_{-n_i - \frac{1}{2}}(x) q^{n_i}
\end{align*} 
where the second equality is by Lemma \ref{lem:computing_M_tilde}. Therefore, the link evaluation $\langle \tilde{\beta} \rangle_{D(H), \mu_R^D}$ equals $q^{\sum_{i=1}^g n_i}$ times the tensor contraction as shown in Figure \ref{fig:M_assignment2} with the latter one being exactly $\Kup{X,b;H}$ described in Section \ref{subsec:computing Kuperberg}. 

Finally, note that $\mu_R^D(v^D) = \tau^{-5}$  and $\mu_R^{D}\left((v^D)^{-1}\right) = \tau$. Choose $\tau^{-3}$ as the square root of $\mu_R^D(v^D) / \mu_R^{D}\left((v^D)^{-1}\right)$.  Hence, the $\hkr$ invariant (the non-refined version) is given by:
\begin{equation}
\HKR{X;D(H), \tau^{-3}} = \tau^{3\sign(\tilde{\beta}) + 2\sum\limits_{i=1}^g w_i} \Kup{X,b;H}. 
\end{equation}
Now choose $\tau^{-1}$ as the sixth root of $\mu_R^D(v^D) / \mu_R^{D}\left((v^D)^{-1}\right)$. Let $\phi$ be the $2$-framing of $X$ such that $\sigma(\phi) = 3\sign(\tilde{\beta}) + 2\sum\limits_{i=1}^g w_i$. Then we get 
\begin{align}
\HKR{X,\phi;D(H), \tau^{-1}} = \Kup{X,b;H}.
\end{align}
\section{Main Results II}
\label{sec:main2}

In this section, the Hopf algebra $H$ is assumed to be factorizable and ribbon. It follows that $H$ is unimodular. We turn to another relation between $Z_{\text{Kup}}$ and $Z_{\text{HKR}}$. It can be viewed as the dual of the relation in 
Theorem \ref{thm:main1}. That is, instead of taking the double of $H$, we take the double $D(X) = X\#\overline{X}$ of the 3-manifold $X$ in $Z_{\text{HKR}}$, where $\overline{X}$ is the manifold $X$ with opposite orientation. 

\begin{theorem}[$=$ Theorem \ref{thm:main22} ]
\label{thm:main2}
Let $H$ be a finite dimensional factorizable ribbon Hopf algebra and $X$ be a closed oriented $3$-manifold, then there exists a framing $b$ of $X$ such that
\begin{align}
\Kup{X,b;H}=\HKR{X\#\overline{X};H}.
\end{align}
\end{theorem}

\begin{figure}
	\centering
	\includegraphics[scale=1]{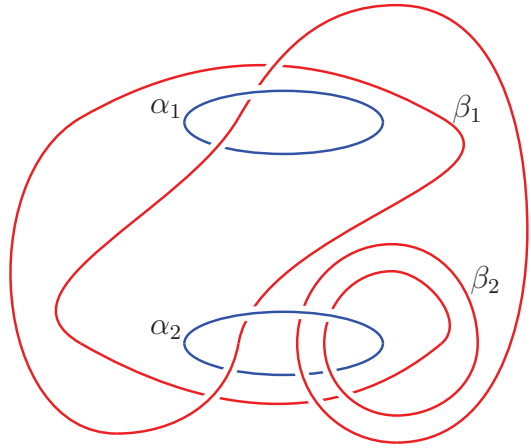}
	\caption{The chain-mail link corresponding to the Heegaard diagram in Figure \ref{fig:special_diagram}}\label{fig:chain-mail}
\end{figure}

The main tool in topology to establish Theorem \ref{thm:main2} is the chain-mail link. A surgery diagram of $X\#\overline{X}$ is obtained from a Heegaard diagram of $X$ by pushing the upper circles into the lower handle body slightly. Then the upper circles and the lower circles form a link $L_{D(X)}$, called a chain-mail link \cite{Roberts}. All these curves are framed by thickening them into thin bands parallel to the Heegaard surface. The framed link $L_{D(X)}$ is a surgery link for $D(X)$.  For instance, Figure \ref{fig:chain-mail} shows the diagram of the chain-mail link for the Heegaard diagram in Figure \ref{fig:special_diagram}.


Note that the signature $\sigma(L_{D(X)})$ of the chain-mail link is always zero \cite{Roberts} and it is possible to choose $\mu_R$ such that $\mu_R(v)\mu_R(v^{-1})= 1$ in
a factorizable ribbon Hopf algebra \cite{Cohen-Westreich}. Hence with such a choice of $\mu_R$ and a suitable choice $\omega(v)$ of a square root of $\mu_R(v)/\mu_R(v^{-1})$, the normalization factor in defining $Z_{\text{HKR}}$ is  
$$[\mu_R(v)\mu_R(v^{-1})]^{-\frac{c(L_{D(X)})}{2}}[\mu_R(v)/\mu_R(v^{-1})]^{-\frac{\sigma(L_{D(X)})}{2}}=1.$$
Thus  $\HKR{X\#\overline{X};H}=\langle L_{D(X)} \rangle_{H,\mu_R}$.

\begin{figure}
\centering
\includegraphics[scale=1]{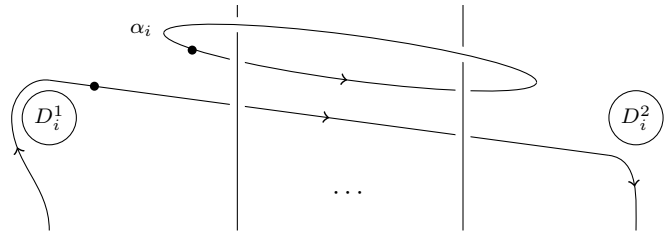}
\caption{Diagram of a chain-main link.} \label{fig:chain_mail_perturb}
\end{figure}

Take $X$ to be $X(\Sigma_g, \beta)$ and choose the framing $b$ to be the one defined in Section \ref{subsec:computing Kuperberg}. We prove $\Kup{X,b;H}=\langle L_{D(X)} \rangle_{H,\mu_R}$. Similar to the proof of Theorem \ref{thm:main1} in Section \ref{subsec:computing hkr}, we perturb the diagram of $L_{D(X)}$, and choose orientation and base point for each component as shown in Figure \ref{fig:chain_mail_perturb}. The following lemma is proved in \cite{Chang-Wang}. Note that we have an extra $S$ factor ($\text{RHS}$ of Figure \ref{fig:resolve_crossing}) compared to the statement in \cite{Chang-Wang}. This is due to the use of a slightly different but equivalent convention in current paper. It is also not hard to verify the lemma directly. 

\begin{lemma}\label{Henningslemma-1}
The equality in Figure \ref{fig:resolve_crossing} holds, where the equality means when the diagram on the $\text{LHS}$ is assigned tensors according to the rules defining $Z_{\text{HKR}}$, then contracting the tensors results in the one on the $\text{RHS}$.
\begin{figure}
\centering
\includegraphics[scale=1]{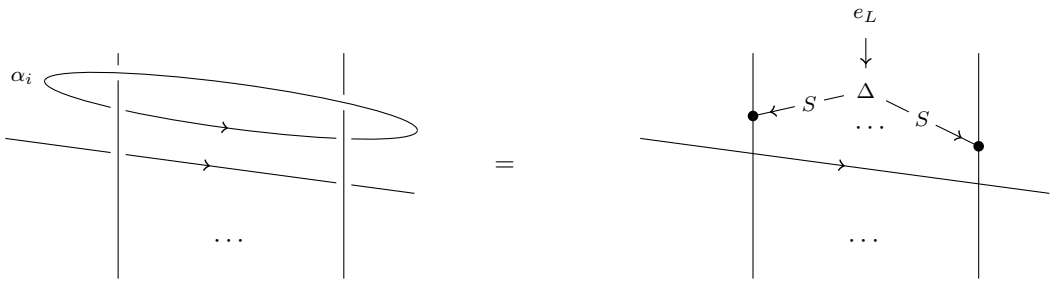}
\caption{Tensor replacement in computing $Z_{\text{HKR}}$. Please be warned that the outgoing legs of the $\Delta$ tensor here is enumerated counterclockwise, in contrast to the default clockwise ordering.}\label{fig:resolve_crossing} 
\end{figure} 
\end{lemma}

Since $H$ is unimodular, we have $e_{n-1/2} = e_L$ for any integer $n$. Lemma \ref{Henningslemma-1} shows that the linking between the lower and upper circles results in the $\Delta$ tensor (Figure \ref{fig:three_tensors}) with an additional $S$ action on each outgoing leg. This effect is the same as assigning the $\Delta$ tensor to the lower circle (with an additional $S$-action). Now for the dot (Figure \ref{fig:resolve_crossing}) corresponding to a crossing $p$, the $S$ powers assigned to it is $S^{-w_p}$. Combining the extra $S$ factor from the previous step, we get $S^{1-w_p}$, which is the correct tensor assigned to the crossing $p$ in the Kuperbeg invariant (see the end of Section \ref{subsec:computing Kuperberg}). Finally,  the $\widetilde{M}$-tensor in the $Z_{\text{HKR}}$ is equal to the $M$-tensor in the $Z_{\text{Kup}}$:
\begin{center}
\includegraphics[scale=1]{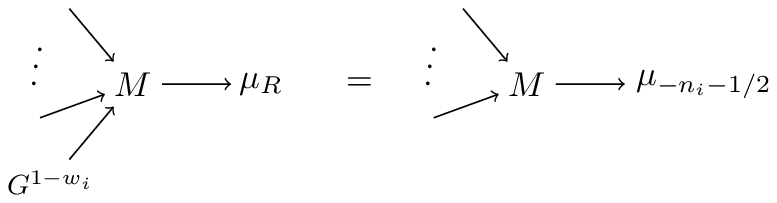}
\end{center}
We get $\langle L_{D(X)} \rangle_{H,\mu_R}=\Kup{X,b;H}$.


\vspace{1cm}
\noindent\textbf{Acknowledgment} $\quad$ The authors would like to thank Greg Kuperberg, Siu-Hung Ng, and Zhenghan Wang for helpful discussions. LC is supported by NSFC Grant No. 11701293. SXC acknowledges the support from the Simons Foundation. 

\bibliographystyle{plain}
\bibliography{KHbib}

\end{document}

%% file: Macros.tex
\theoremstyle{definition}
\newtheorem{theorem}{Theorem}[section]

\newtheorem{lemma}[theorem]{Lemma}

\newtheorem{remark}[theorem]{Remark}

\def\Real{\mathbb{R}}
\def\Integer{\mathbb{Z}}
\def\Complex{\mathbb{C}}
\def\CatC{\mathcal{C}}

\def\CatT{\mathcal{T}}

\def\CatG{\mathcal{G}}

\def\TQFT{\mathrm{TQFT}}

\def\Sphere{\mathbb{S}}

\def\SO{\mathrm{SO}}
\def\Spin{\mathrm{Spin}}
\def\Dop{\Delta^{{\footnotesize \text{op}}}}
\def\cop{{\footnotesize \text{cop}}}
\def\qt{\text{quasitriangular}}
\def\Drinfeld{\text{Drinfeld}}
\def\db{\text{double balanced}}

\def\hkr{\text{Hennings-Kauffman-Radford}}
\def\wrt{\text{Witten-Reshetikhin-Turaev}}

\def\HKR#1{Z_{\text{HKR}}(#1)}
\def\Kup#1{Z_{\text{Kup}}(#1)}
\def\sign{\text{sign}}